\documentstyle[12pt,fleqn]{article}

\begin{document}

\title{Complex Numbers in 6 Dimensions}

\author{Silviu Olariu
\thanks{e-mail: olariu@ifin.nipne.ro}\\
Institute of Physics and Nuclear Engineering, Tandem Laboratory\\
76900 Magurele, P.O. Box MG-6, Bucharest, Romania}

\date{4 August 2000}

\maketitle

\abstract
Two distinct systems of commutative complex numbers in 6 dimensions of the
polar and planar types having the form
$u=x_0+h_1x_1+h_2x_2+h_3x_3+h_4x_4+h_5x_5$, are described in this work, where
the variables $x_0, x_1, x_2, x_3, x_4, x_5$ are real numbers.  The polar
6-complex numbers introduced in this paper can be specified by the modulus $d$,
the amplitude $\rho$, and the polar angles $\theta_+,
\theta_-$, the planar angle $\psi_1$, and the azimuthal angles $\phi_1,
\phi_2$. The 
planar 6-complex numbers introduced in this paper can be specified by the
modulus $d$, the amplitude $\rho$, the planar angles $\psi_1, \psi_2$, and the
azimuthal angles $\phi_1, \phi_2, \phi_3$.  Exponential and trigonometric forms
are given for the 6-complex numbers.  The 6-complex functions defined by series
of powers are analytic, and the partial derivatives of the components of the
6-complex functions are closely related.  The integrals of polar 6-complex
functions are independent of path in regions where the functions are regular.
The fact that the exponential form of ther 6-complex numbers depends on cyclic
variables leads to the concept of pole and residue for integrals on closed
paths. The polynomials of polar 6-complex variables can be written as products
of linear or quadratic factors, the polynomials of planar 6-complex variables
can always be written as products of linear factors, although the factorization
is not unique.

\endabstract

\section{Introduction}

A regular, two-dimensional complex number $x+iy$ 
can be represented geometrically by the modulus $\rho=(x^2+y^2)^{1/2}$ and 
by the polar angle $\theta=\arctan(y/x)$. The modulus $\rho$ is multiplicative
and the polar angle $\theta$ is additive upon the multiplication of ordinary 
complex numbers.

The quaternions of Hamilton are a system of hypercomplex numbers
defined in four dimensions, the
multiplication being a noncommutative operation, \cite{1} 
and many other hypercomplex systems are
possible, \cite{2a}-\cite{2b} but these hypercomplex systems 
do not have all the required properties of regular, 
two-dimensional complex numbers which rendered possible the development of the 
theory of functions of a complex variable.

Two distinct systems of commutative complex numbers in 6 dimensions 
having the form $u=x_0+h_1x_1+h_2x_2+h_3x_3+h_4x_4+h_5x_5$ are
described in this work, 
for which the multiplication is associative and
commutative,
where the variables $x_0, x_1, x_2, x_3, x_4, x_5$ are
real numbers.  
The first type of 6-complex numbers described in this article is
characterized by the presence
of two polar axes, so that
these numbers will be called polar 6-complex numbers. 
The other type of 6-complex numbers described in this paper 
will be called planar n-complex numbers. 

The polar 6-complex numbers introduced in this paper can be specified by the
modulus $d$, the amplitude $\rho$, and the polar angles $\theta_+,
\theta_-$, the planar angle $\psi_1$, and the azimuthal angles $\phi_1,
\phi_2$. The 
planar 6-complex numbers introduced in this paper can be specified by the
modulus $d$, the amplitude $\rho$, the planar angles $\psi_1, \psi_2$, and the
azimuthal angles $\phi_1, \phi_2, \phi_3$.  Exponential and trigonometric forms
are given for the 6-complex numbers.  The 6-complex functions defined by series
of powers are analytic, and the partial derivatives of the components of the
6-complex functions are closely related.  The integrals of polar 6-complex
functions are independent of path in regions where the functions are regular.
The fact that the exponential form of ther 6-complex numbers depends on cyclic
variables leads to the concept of pole and residue for integrals on closed
paths. The polynomials of polar 6-complex variables can be written as products
of linear or quadratic factors, the polynomials of planar 6-complex variables
can always be written as products of linear factors, although the factorization
is not unique.

This paper belongs to a series of studies on commutative complex numbers in $n$
dimensions. \cite{2c}
The polar 6-complex numbers described in this paper are a particular case for
$n=6$ of the polar hypercomplex numbers in $n$ dimensions, and the
planar 6-complex numbers described in this section are a particular case for
$n=6$ of the planar hypercomplex numbers in $n$ dimensions.\cite{2c},\cite{2d}

\section{Polar complex numbers in 6 dimensions}

\subsection{Operations with polar complex numbers in $6$ dimensions}

The polar hypercomplex number $u$ in 6 dimensions 
is represented as  
\begin{equation}
u=x_0+h_1x_1+h_2x_2+h_3x_3+h_4x_4+h_5x_5. 
\label{6ch1a}
\end{equation}
The multiplication rules for the bases 
$h_1, h_2, h_3, h_4, h_5 $ are 
\begin{eqnarray}
\lefteqn{h_1^2=h_2,\;h_2^2=h_4,\;h_3^2=1,\;h_4^2=h_2,\;h_5^2=h_4,\;
h_1h_2=h_3,\;h_1h_3=h_4,\;h_1h_4=h_5,\nonumber}\\
&&\;h_1h_5=1,\;h_2h_3=h_5,\;h_2h_4=1,\;h_2h_5=h_1,\;h_3h_4=h_1,
\;h_3h_5=h_2,\;h_4h_5=h_3.
\label{6ch1}
\end{eqnarray}
The significance of the composition laws in Eq.  (\ref{6ch1}) can be understood
by representing the bases $h_j, h_k$ by points on a circle at the angles
$\alpha_j=\pi j/3,\alpha_k=\pi k/3$, as shown in Fig. 1, and the product $h_j
h_k$ by the point of the circle at the angle $\pi (j+k)/3$. If $2\pi\leq \pi
(j+k)/3<4\pi$, the point represents the basis $h_l$ of angle
$\alpha_l=\pi(j+k)/3-2\pi$.

The sum of the 6-complex numbers $u$ and $u^\prime$ is
\begin{equation}
u+u^\prime=x_0+x^\prime_0+h_1(x_1+x^\prime_1)+h_1(x_2+x^\prime_2)
+h_3(x_3+x^\prime_3)+h_4(x_4+x^\prime_4)+h_5(x_5+x^\prime_5).
\label{6ch2}
\end{equation}
The product of the numbers $u, u^\prime$ is
\begin{equation}
\begin{array}{l}
uu^\prime=x_0 x_0^\prime +x_1x_5^\prime+x_2 x_4^\prime+x_3x_3^\prime
+x_4x_2^\prime+x_5 x_1^\prime\\
+h_1(x_0 x_1^\prime+x_1x_0^\prime+x_2x_5^\prime+x_3x_4^\prime
+x_4 x_3^\prime+x_5 x_2^\prime) \\
+h_2(x_0 x_2^\prime+x_1x_1^\prime+x_2x_0^\prime+x_3x_5^\prime
+x_4 x_4^\prime+x_5 x_3^\prime) \\
+h_3(x_0 x_3^\prime+x_1x_2^\prime+x_2x_1^\prime+x_3x_0^\prime
+x_4 x_5^\prime+x_5 x_4^\prime) \\
+h_4(x_0 x_4^\prime+x_1x_3^\prime+x_2x_2^\prime+x_3x_1^\prime
+x_4 x_0^\prime+x_5 x_5^\prime) \\
+h_5(x_0 x_5^\prime+x_1x_4^\prime+x_2x_3^\prime+x_3x_2^\prime
+x_4 x_1^\prime+x_5 x_0^\prime).
\end{array}
\label{6ch3}
\end{equation}

The relation between the variables $v_+,v_-,v_1,\tilde v_1, v_2, \tilde v_2$
and $x_0,x_1,x_2,x_3,x_4,x_5$ are
\begin{equation}
\left(
\begin{array}{c}
v_+\\
v_-\\
v_1\\
\tilde v_1\\
v_2\\
\tilde v_2\\
\end{array}\right)
=\left(
\begin{array}{cccccc}
1&1&1&1&1&1\\
1&-1&1&-1&1&-1\\
1&\frac{1}{2} &-\frac{1}{2} &-1 &-\frac{1}{2} &\frac{1}{2} \\
0& \frac{\sqrt{3}}{2}&\frac{\sqrt{3}}{2} &0 &-\frac{\sqrt{3}}{2}
&-\frac{\sqrt{3}}{2} \\ 
1& -\frac{1}{2}&-\frac{1}{2} &1 &-\frac{1}{2} &-\frac{1}{2} \\
0&\frac{\sqrt{3}}{2} &-\frac{\sqrt{3}}{2} &0 &\frac{\sqrt{3}}{2}
&-\frac{\sqrt{3}}{2} \\ 
\end{array}
\right)
\left(
\begin{array}{c}
x_0\\
x_1\\
x_2\\
x_3\\
x_4\\
x_5
\end{array}
\right).
\label{6ch9e}
\end{equation}
The other variables are $v_4=v_2, \tilde v_4=-\tilde v_2,
v_5=v_1, \tilde v_5=-\tilde v_1$. 
The variables $v_+, v_-, v_1, \tilde v_1, v_2, \tilde v_2$ will be called
canonical polar 6-complex variables.

\subsection{Geometric representation of polar complex numbers in 6 dimensions}

The 6-complex number $u=x_0+h_1x_1+h_2x_2+h_3x_3+h_4x_4+h_5x_5$
is represented by 
the point $A$ of coordinates $(x_0,x_1,x_2,x_3,x_4,x_5)$. 
The distance from the origin $O$ of the 6-dimensional space to the point $A$
has the expression 
\begin{equation}
d^2=x_0^2+x_1^2+x_2^2+x_3^2+x_4^2+x_5^2.
\label{6ch10}
\end{equation}
The distance $d$ is called modulus of the 6-complex number 
$u$, and is designated by $d=|u|$.
The modulus has the property that
\begin{equation}
|u^\prime u^{\prime\prime}|\leq \sqrt{6}|u^\prime||u^{\prime\prime}| .
\label{6ch79}
\end{equation}

The exponential and trigonometric forms of the 6-complex number $u$ can be
obtained conveniently in a rotated system of axes defined by a transformation
which has the form
\begin{equation}
\left(
\begin{array}{c}
\xi_+\\
\xi_-\\
\xi_1\\
\tilde \xi_1\\
\xi_2\\
\tilde \xi_2\\
\end{array}\right)
=\left(
\begin{array}{cccccc}
\frac{1}{\sqrt{6}}&\frac{1}{\sqrt{6}}&
\frac{1}{\sqrt{6}}&\frac{1}{\sqrt{6}}&\frac{1}{\sqrt{6}}&\frac{1}{\sqrt{6}}\\
\frac{1}{\sqrt{6}}&-\frac{1}{\sqrt{6}}&\frac{1}{\sqrt{6}}&-\frac{1}{\sqrt{6}}&
\frac{1}{\sqrt{6}}&-\frac{1}{\sqrt{6}}\\
\frac{\sqrt{3}}{3}&\frac{\sqrt{3}}{6} &
-\frac{\sqrt{3}}{6} &-\frac{\sqrt{3}}{3} &
-\frac{\sqrt{3}}{6} &\frac{\sqrt{3}}{6} \\
0& \frac{1}{2}&\frac{1}{2} &0 &-\frac{1}{2} &-\frac{1}{2} \\
\frac{\sqrt{3}}{3}& -\frac{\sqrt{3}}{6}&
-\frac{\sqrt{3}}{6} &\frac{\sqrt{3}}{3} &
-\frac{\sqrt{3}}{6} &-\frac{\sqrt{3}}{6} \\
0&\frac{1}{2} &-\frac{1}{2} &0 &\frac{1}{2} &-\frac{1}{2} \\
\end{array}
\right)
\left(
\begin{array}{c}
x_0\\
x_1\\
x_2\\
x_3\\
x_4\\
x_5
\end{array}
\right).
\label{6ch9ee}
\end{equation}
The lines of the matrices in Eq. (\ref{6ch9ee})  gives the components
of the 6 basis vectors of the new system of axes. These vectors have unit
length and are orthogonal to each other.
The relations between the two sets of variables are
\begin{equation}
v_+= \sqrt{6}\xi_+ ,  v_-=  \sqrt{6}\xi_-,  
v_k= \sqrt{3}\xi_k , \tilde v_k= \sqrt{3}\eta_k, k=1,2 .
\label{6ch12b}
\end{equation}

The radius $\rho_k$ and the azimuthal angle $\phi_k$ in the plane of the axes
$v_k,\tilde v_k$ are
\begin{equation}
\rho_k^2=v_k^2+\tilde v_k^2, \:\cos\phi_k=v_k/\rho_k,
\:\sin\phi_k=\tilde v_k/\rho_k, 0\leq \phi_k<2\pi ,\;k=1,2,
\label{6ch19a}
\end{equation}
so that there are 2 azimuthal angles.
The planar angle $\psi_1$ is
\begin{equation}
\tan\psi_1=\rho_1/\rho_2, 0\leq\psi_1\leq\pi/2.
\label{6ch19b}
\end{equation}
There is a polar angle $\theta_+$, 
\begin{equation}
\tan\theta_+=\frac{\sqrt{2}\rho_1}{v_+}, 0\leq\theta_+\leq\pi , 
\label{6ch19c}
\end{equation}
and there is also a polar angle $\theta_-$,
\begin{equation}
\tan\theta_-=\frac{\sqrt{2}\rho_1}{v_-}, 0\leq\theta_-\leq\pi .
\label{6ch19d}
\end{equation}
The amplitude of a 6-complex number $u$ is
\begin{equation}
\rho=\left(v_+v_-\rho_1^2 \rho_2^2\right)^{1/6}.
\label{6ch50aa}
\end{equation}
It can be checked that
\begin{equation}
d^2=\frac{1}{6}v_+^2+\frac{1}{6}v_-^2
+\frac{1}{3}(\rho_1^2+\rho_2^2) .
\label{6ch17}
\end{equation}

If $u=u^\prime u^{\prime\prime}$, the parameters 
of the hypercomplex numbers are related by
\begin{equation}
v_+=v_+^\prime v_+^{\prime\prime},  
\label{6ch21a}
\end{equation}
\begin{equation}
\tan\theta_+=\frac{1}{\sqrt{2}}\tan\theta_+^\prime \tan\theta_+^{\prime\prime},
\label{6ch21c}
\end{equation}
\begin{equation}
v_-=v_-^\prime v_-^{\prime\prime},  
\label{6ch21f}
\end{equation}
\begin{equation}
\tan\theta_-=\frac{1}{\sqrt{2}}\tan\theta_-^\prime \tan\theta_-^{\prime\prime},
\label{6ch21g}
\end{equation}
\begin{equation}
\tan\psi_1=\tan\psi_1^\prime \tan\psi_1^{\prime\prime},
\label{6ch21d}
\end{equation}
\begin{equation}
\rho_k=\rho_k^\prime\rho_k^{\prime\prime}, 
\label{6ch21b}
\end{equation}
\begin{equation}
\phi_k=\phi_k^\prime+\phi_k^{\prime\prime},  
\label{6ch21e}
\end{equation}
\begin{equation}
v_k=v_k^\prime v_k^{\prime\prime}-\tilde v_k^\prime \tilde
v_k^{\prime\prime},\;
\tilde v_k=v_k^\prime 
\tilde v_k^{\prime\prime}+\tilde v_k^\prime v_k^{\prime\prime},
\label{6ch22}
\end{equation}
\begin{equation}
\rho=\rho^\prime\rho^{\prime\prime} ,
\label{6ch24}
\end{equation}
where $k=1,2$.

The 6-complex number
$u=x_0+h_1x_1+h_2x_2+h_3x_3+h_4x_4+h_5x_5$ can be  represented by the matrix
\begin{equation}
U=\left(
\begin{array}{cccccc}
x_0     &   x_1     &   x_2   &   x_3  &  x_4  & x_5\\
x_5     &   x_0     &   x_1   &   x_2  &  x_3  & x_4\\
x_4     &   x_5     &   x_0   &   x_1  &  x_2  & x_3\\
x_3     &   x_4     &   x_5   &   x_0  &  x_1  & x_2\\
x_2     &   x_3     &   x_4   &   x_5  &  x_0  & x_1\\
x_1     &   x_2     &   x_3   &   x_4  &  x_5  & x_0\\
\end{array}
\right).
\label{6ch24b}
\end{equation}
The product $u=u^\prime u^{\prime\prime}$ is
represented by the matrix multiplication $U=U^\prime U^{\prime\prime}$.

\subsection{The polar 6-dimensional cosexponential functions}

The polar cosexponential functions in 6 dimensions are
\begin{equation}
g_{6k}(y)=\sum_{p=0}^\infty y^{k+6p}/(k+6p)!, 
\label{6ch29}
\end{equation}
for $ k=0,...,5$.
The polar cosexponential functions $g_{6k}$ of even index $k$ are
even functions, $g_{6,2p}(-y)=g_{6,2p}(y)$, 
and the polar cosexponential functions of odd index $k$
are odd functions, $g_{6,2p+1}(-y)=-g_{6,2p+1}(y)$, $p=0,1,2$. 

It can be checked that
\begin{equation}
\sum_{k=0}^{5}g_{6k}(y)=e^y,
\label{6ch29a}
\end{equation}
\begin{equation}
\sum_{k=0}^{5}(-1)^k g_{6k}(y)=e^{-y}.
\label{6ch29b}
\end{equation}
The exponential function of the quantity $h_k y$ is
\begin{equation}
\begin{array}{l}
e^{h_1 y}=g_{60}(y)+h_1g_{61}(y)+h_2g_{62}(y)
+h_3g_{63}(y)+h_4g_{64}(y)+h_5g_{65}(y),\\
e^{h_2 y}=g_{60}(y)+g_{63}(y)+h_2\{g_{61}(y)
+g_{64}(y)\}+h_4\{g_{62}(y)+g_{65}(y)\},\\
e^{h_3 y}=g_{60}(y)+g_{62}(y)+g_{64}(y)+h_3\{g_{61}(y)+g_{63}(y)+g_{65}(y)\},\\
e^{h_4 y}=g_{60}(y)+g_{63}(y)+h_2\{g_{62}(y)
+g_{65}(y)\}+h_4\{g_{61}(y)+g_{64}(y)\},\\
e^{h_5 y}=g_{60}(y)+h_1g_{65}(y)+h_2g_{64}(y)
+h_3g_{63}(y)+h_4g_{62}(y)+h_5g_{61}(y).\\
\end{array}
\label{6ch28b}
\end{equation}
The relations for $h_2$ and $h_4$ can be written equivalently as 
$e^{h_2 y}=g_{30}+h_2g_{31}+h_4 g_{32}, e^{h_4 y}=g_{30}+h_2g_{32}+h_4 g_{31}$,
and the relation for $h_3$ can be written as $e^{h_3 y}=g_{20}+h_3g_{21}$,
which is the same as $e^{h_3 y}=\cosh y+h_3\sinh y$.

The expressions of the polar 6-dimensional cosexponential functions are
\begin{equation}
\begin{array}{l}
g_{60}(y)=\frac{1}{3}\cosh y +\frac{2}{3}\cosh\frac{y}{2}
\cos\frac{\sqrt{3}}{2}y,\\
g_{61}(y)=\frac{1}{3}\sinh y +\frac{1}{3}\sinh\frac{y}{2}
\cos\frac{\sqrt{3}}{2}y
+\frac{\sqrt{3}}{3}\cosh\frac{y}{2}\sin\frac{\sqrt{3}}{2}y,\\
g_{62}(y)=\frac{1}{3}\cosh y -\frac{1}{3}\cosh\frac{y}{2}
\cos\frac{\sqrt{3}}{2}y
+\frac{\sqrt{3}}{3}\sinh\frac{y}{2}\sin\frac{\sqrt{3}}{2}y,\\
g_{63}(y)=\frac{1}{3}\sinh y -\frac{2}{3}\sinh\frac{y}{2}
\cos\frac{\sqrt{3}}{2}y,\\
g_{64}(y)=\frac{1}{3}\cosh y -\frac{1}{3}\cosh\frac{y}{2}
\cos\frac{\sqrt{3}}{2}y
-\frac{\sqrt{3}}{3}\sinh\frac{y}{2}\sin\frac{\sqrt{3}}{2}y,\\
g_{65}(y)=\frac{1}{3}\sinh y +\frac{1}{3}\sinh\frac{y}{2}
\cos\frac{\sqrt{3}}{2}y
-\frac{\sqrt{3}}{3}\cosh\frac{y}{2}\sin\frac{\sqrt{3}}{2}y.
\end{array}
\label{6ch30}
\end{equation}
The cosexponential functions (\ref{6ch30}) can be written as
\begin{equation}
g_{6k}(y)=\frac{1}{6}\sum_{l=0}^{5}\exp\left[y\cos\left(\frac{2\pi l}{6}\right)
\right]
\cos\left[y\sin\left(\frac{2\pi l}{6}\right)-\frac{2\pi kl}{6}\right], 
\label{6ch30x}
\end{equation}
for $k=0,...,5$.
The graphs of the polar 6-dimensional cosexponential functions are shown in Fig
2. 

It can be checked that
\begin{equation}
\sum_{k=0}^{5}g_{6k}^2(y)=\frac{1}{3}\cosh 2y +\frac{2}{3}\cosh y.
\label{6ch34a}
\end{equation}

The addition theorems for the polar 6-dimensional cosexponential functions are
\begin{eqnarray}
\lefteqn{\begin{array}{l}
g_{60}(y+z)=g_{60}(y)g_{60}(z)+g_{61}(y)g_{65}(z)+g_{62}(y)g_{64}(z)
+g_{63}(y)g_{63}(z)+g_{64}(y)g_{62}(z)+g_{65}(y)g_{61}(z) ,\\
g_{61}(y+z)=g_{60}(y)g_{61}(z)+g_{61}(y)g_{60}(z)+g_{62}(y)g_{65}(z)
+g_{63}(y)g_{64}(z)+g_{64}(y)g_{63}(z)+g_{65}(y)g_{62}(z) ,\\
g_{62}(y+z)=g_{60}(y)g_{62}(z)+g_{61}(y)g_{61}(z)+g_{62}(y)g_{60}(z)
+g_{63}(y)g_{65}(z)+g_{64}(y)g_{64}(z)+g_{65}(y)g_{63}(z) ,\\
g_{63}(y+z)=g_{60}(y)g_{63}(z)+g_{61}(y)g_{62}(z)+g_{62}(y)g_{61}(z)
+g_{63}(y)g_{60}(z)+g_{64}(y)g_{65}(z)+g_{65}(y)g_{64}(z) ,\\
g_{64}(y+z)=g_{60}(y)g_{64}(z)+g_{61}(y)g_{63}(z)+g_{62}(y)g_{62}(z)
+g_{63}(y)g_{61}(z)+g_{64}(y)g_{60}(z)+g_{65}(y)g_{65}(z) ,\\
g_{65}(y+z)=g_{60}(y)g_{65}(z)+g_{61}(y)g_{64}(z)+g_{62}(y)g_{63}(z)
+g_{63}(y)g_{62}(z)+g_{64}(y)g_{61}(z)+g_{65}(y)g_{60}(z) .
\end{array}\nonumber}\\
&&
\label{6ch35a}
\end{eqnarray}
It can be shown that
\begin{eqnarray}
\begin{array}{l}
\{g_{60}(y)+h_1g_{61}(y)+h_2g_{62}(y)+h_3g_{63}(y)
+h_4g_{64}(y)+h_5g_{65}(y)\}^l\\
\hspace{0.5cm}=g_{60}(ly)+h_1g_{61}(ly)+h_2g_{62}(ly)
+h_3g_{63}(ly)+h_4g_{64}(ly)+h_5g_{65}(ly),\\
\{g_{60}(y)+g_{63}(y)+h_2\{g_{61}(y)+g_{64}(y)\}
+h_4\{g_{62}(y)+g_{65}(y)\}\}^l\\
\hspace{0.5cm}=g_{60}(ly)+g_{63}(ly)+h_2\{g_{61}(ly)
+g_{64}(ly)\}+h_4\{g_{62}(ly)+g_{65}(ly)\},\\
\{g_{60}(y)+g_{62}(y)+g_{64}(y)+h_3\{g_{61}(y)+g_{63}(y)+g_{65}(y)\}\}^l\\
\hspace{0.5cm}=g_{60}(ly)+g_{62}(ly)+g_{64}(ly)
+h_3\{g_{61}(ly)+g_{63}(ly)+g_{65}(ly)\},\\
\{g_{60}(y)+g_{63}(y)+h_2\{g_{62}(y)+g_{65}(y)\}
+h_4\{g_{61}(y)+g_{64}(y)\}\}^l\\
\hspace{0.5cm}=g_{60}(ly)+g_{63}(ly)+h_2\{g_{62}(ly)
+g_{65}(ly)\}+h_4\{g_{61}(ly)+g_{64}(ly)\},\\
\{g_{60}(y)+h_1g_{65}(y)+h_2g_{64}(y)+h_3g_{63}(y)
+h_4g_{62}(y)+h_5g_{61}(y)\}^l\\
\hspace{0.5cm}=g_{60}(ly)+h_1g_{65}(ly)
+h_2g_{64}(ly)+h_3g_{63}(ly)+h_4g_{62}(ly)+h_5g_{61}(ly).\\
\end{array}
\label{6ch37b}
\end{eqnarray}

The derivatives of the polar cosexponential functions
are related by
\begin{equation}
\frac{dg_{60}}{du}=g_{65}, \:
\frac{dg_{61}}{du}=g_{60}, \:
\frac{dg_{62}}{du}=g_{61}, \:
\frac{dg_{63}}{du}=g_{62} ,\:
\frac{dg_{64}}{du}=g_{63}, \:
\frac{dg_{65}}{du}=g_{64} .
\label{6ch45}
\end{equation}

\subsection{Exponential and trigonometric forms of polar 6-complex numbers}

The exponential and trigonometric forms of polar 6-complex
numbers can be expressed with the aid of the hypercomplex bases 
\begin{equation}
\left(
\begin{array}{c}
e_+\\
e_-\\
e_1\\
\tilde e_1\\
e_2\\
\tilde e_2\\
\end{array}\right)
=\left(
\begin{array}{cccccc}
\frac{1}{6}&\frac{1}{6}&\frac{1}{6}&\frac{1}{6}&\frac{1}{6}&\frac{1}{6}\\
\frac{1}{6}&-\frac{1}{6}&\frac{1}{6}&-\frac{1}{6}&\frac{1}{6}&-\frac{1}{6}\\
\frac{1}{3}&\frac{1}{6} &-\frac{1}{6} &
-\frac{1}{3} &-\frac{1}{6} &\frac{1}{6} \\
0& \frac{\sqrt{3}}{6}&\frac{\sqrt{3}}{6} &
0 &-\frac{\sqrt{3}}{6} &-\frac{\sqrt{3}}{6} \\
\frac{1}{3}& -\frac{1}{6}&-\frac{1}{6} &
\frac{1}{3} &-\frac{1}{6} &-\frac{1}{6} \\
0&\frac{\sqrt{3}}{6} &-\frac{\sqrt{3}}{6} &
0 &\frac{\sqrt{3}}{6} &-\frac{\sqrt{3}}{6} \\
\end{array}
\right)
\left(
\begin{array}{c}
1\\
h_1\\
h_2\\
h_3\\
h_4\\
h_5
\end{array}
\right).
\label{6che11}
\end{equation}

The multiplication relations for these bases are
\begin{eqnarray}
\lefteqn{e_+^2=e_+,\; e_-^2=e_-,\; e_+e_-=0,\; e_+e_k=0,\; e_+\tilde e_k=0,\;
e_-e_k=0,\; 
e_-\tilde e_k=0,\nonumber}\\ 
&&e_k^2=e_k,\; \tilde e_k^2=-e_k,\; e_k \tilde e_k=\tilde e_k ,\; e_ke_l=0,\;
e_k\tilde e_l=0,\; \tilde e_k\tilde e_l=0,\; k,l=1,2, k\not=l.\; 
\label{6che12a}
\end{eqnarray}
The bases have the property that
\begin{equation}
e_+ +e_- + e_1+e_2 =1.
\label{6ch47b}
\end{equation}
The moduli of the new bases are
\begin{equation}
|e_+|=\frac{1}{\sqrt{6}},\; |e_-|=\frac{1}{\sqrt{6}},\; 
|e_k|=\frac{1}{\sqrt{3}},\; |\tilde e_k|=\frac{1}{\sqrt{3}}, k=1,2.
\label{6che12c}
\end{equation}

It can be shown that
\begin{eqnarray}
\lefteqn{x_0+h_1x_1+h_2x_2+h_3x_3+h_4x_4+h_5x_5\nonumber}\\
&&=e_+v_+ + e_-v_- +e_1 v_1+\tilde e_1 \tilde v_1
+e_2 v_2+\tilde e_2 \tilde v_2.
\label{6che13a}
\end{eqnarray}
The ensemble $e_+, e_-, e_1, \tilde e_1, e_2, \tilde e_2$ will be called the
canonical polar 6-complex base, and Eq. (\ref{6che13a}) gives the canonical
form of the polar 6-complex number.

The exponential form of the 6-complex number $u$ is 
\begin{eqnarray}\lefteqn{
u=\rho\exp\left\{\frac{1}{6}(h_1+h_2+h_3+h_4+h_5)
\ln\frac{\sqrt{2}}{\tan\theta_+}
-\frac{1}{6}(h_1-h_2+h_3-h_4+h_5)
\ln\frac{\sqrt{2}}{\tan\theta_-}\right.\nonumber}\\
&&
\left.+\frac{1}{6}(h_1+h_2-2h_3+h_4+h_5)\ln\tan\psi_1
+\tilde e_1\phi_1+\tilde e_2\phi_2 \right\},
\label{6ch50a}
\end{eqnarray}
for $0<\theta_+<\pi/2, 0<\theta_-<\pi/2$.

The
trigonometric form of the 6-complex number $u$ is
\begin{eqnarray}
\lefteqn{u=d
\sqrt{3}
\left(\frac{1}{\tan^2\theta_+}+\frac{1}{\tan^2\theta_-}+1
+\frac{1}{\tan^2\psi_1}\right)^{-1/2}\nonumber}\\
&&\left(\frac{e_+\sqrt{2}}{\tan\theta_+}+\frac{e_-\sqrt{2}}{\tan\theta_-}
+e_1+\frac{e_2}{\tan\psi_1}\right)
\exp\left(\tilde e_1\phi_1+\tilde e_2\phi_2\right).
\label{6ch52a}
\end{eqnarray}

The modulus $d$ and the amplitude $\rho$ are related by
\begin{eqnarray}
\lefteqn{d=\rho \frac{2^{1/3}}{\sqrt{6}}
\left(\tan\theta_+\tan\theta_-
\tan^2\psi_1\right)^{1/6}\nonumber}\\
&&\left(\frac{1}{\tan^2\theta_+}+\frac{1}{\tan^2\theta_-}+1
+\frac{1}{\tan^2\psi_1}\right)^{1/2}.
\label{6ch53a}
\end{eqnarray}

\subsection{Elementary functions of a polar 6-complex variable}

The logarithm and power functions of the 6-complex number $u$ exist for
$v_+>0, v_->0$, which means that $0<\theta_+<\pi/2, 0<\theta_-<\pi/2$,
and are given by
\begin{eqnarray}\lefteqn{
\ln u=\ln \rho+
\frac{1}{6}(h_1+h_2+h_3+h_4+h_5)\ln\frac{\sqrt{2}}{\tan\theta_+}
-\frac{1}{6}(h_1-h_2+h_3-h_4+h_5)\ln\frac{\sqrt{2}}{\tan\theta_-}
\nonumber}\\
&&
+\frac{1}{6}(h_1+h_2-2h_3+h_4+h_5)\ln\tan\psi_1
+\tilde e_1\phi_1+\tilde e_2\phi_2 ,
\label{6ch56a}
\end{eqnarray}
\begin{equation}
u^m=e_+ v_+^m+e_- v_-^m +
\rho_1^m(e_1\cos m\phi_1+\tilde e_1\sin m\phi_1)
+\rho_2^m(e_2\cos m\phi_2+\tilde e_2\sin m\phi_2).
\label{6ch59a}
\end{equation}

The exponential of the 6-complex variable $u$ is 
\begin{eqnarray}
e^u=e_+e^{v_+} + e_-e^{v_-} 
+e^{v_1}\left(e_1 \cos \tilde v_1+\tilde e_1 \sin\tilde v_1\right)
+e^{v_2}\left(e_2 \cos \tilde v_2+\tilde e_2 \sin\tilde v_2\right).
\label{6ch73a}
\end{eqnarray}
The trigonometric functions of the
6-complex variable $u$ are 
\begin{equation}
\cos u=e_+\cos v_+ + e_-\cos v_- 
+\sum_{k=1}^{2}\left(e_k \cos v_k\cosh \tilde v_k
-\tilde e_k \sin v_k\sinh\tilde v_k\right),
\label{6ch74a}
\end{equation}
\begin{equation}
\sin u=e_+\sin v_+ + e_-\sin v_- 
+\sum_{k=1}^{2}\left(e_k \sin v_k\cosh \tilde v_k
+\tilde e_k \cos v_k\sinh\tilde v_k\right).
\label{6ch74b}
\end{equation}

The hyperbolic functions of the
6-complex variable $u$ are
\begin{equation}
\cosh u=e_+\cosh v_+ + e_-\cosh v_- 
+\sum_{k=1}^{2}\left(e_k \cosh v_k\cos \tilde v_k
+\tilde e_k \sinh v_k\sin\tilde v_k\right),
\label{6ch75a}
\end{equation}
\begin{equation}
\sinh u=e_+\sinh v_+ + e_-\sinh v_- 
+\sum_{k=1}^{2}\left(e_k \sinh v_k\cos \tilde v_k
+\tilde e_k \cosh v_k\sin\tilde v_k\right).
\label{6ch75b}
\end{equation}

\subsection{Power series of 6-complex numbers}

A power series of the 6-complex variable $u$ is a series of the form
\begin{equation}
a_0+a_1 u + a_2 u^2+\cdots +a_l u^l+\cdots .
\label{6ch83}
\end{equation}
Since
\begin{equation}
|au^l|\leq 6^{l/2} |a| |u|^l ,
\label{6ch82}
\end{equation}
the series is absolutely convergent for 
\begin{equation}
|u|<c,
\label{6ch86}
\end{equation}
where 
\begin{equation}
c=\lim_{l\rightarrow\infty} \frac{|a_l|}{\sqrt{6}|a_{l+1}|} .
\label{6ch87}
\end{equation}

If $a_l=\sum_{p=0}^{5}h_p a_{lp}$, where $h_0=1$, and
\begin{equation}
A_{l+}=\sum_{p=0}^{5}a_{lp},
\label{6ch88a}
\end{equation}
\begin{equation}
A_{l-}=\sum_{p=0}^{5}(-1)^p a_{lp},
\label{6ch88d}
\end{equation}
\begin{equation}
A_{lk}=\sum_{p=0}^{5}a_{lp}\cos\frac{\pi kp}{3},
\label{6ch88b}
\end{equation}
\begin{equation}
\tilde A_{lk}=\sum_{p=0}^{5}a_{lp}\sin\frac{\pi kp}{3},
\label{6ch88c}
\end{equation}
for $k=1,2$,   
the series (\ref{6ch83}) can be written as
\begin{equation}
\sum_{l=0}^\infty \left[
e_+A_{l+}v_+^l+e_-A_{l-}v_-^l+\sum_{k=1}^{2}
(e_k A_{lk}+\tilde e_k\tilde A_{lk})(e_k v_k+\tilde e_k\tilde v_k)^l 
\right].
\label{6ch89a}
\end{equation}

The series in Eq. (\ref{6ch83}) is absolutely convergent for 
\begin{equation}
|v_+|<c_+,\:
|v_-|<c_-,\:
\rho_k<c_k, k=1,2,
\label{6ch90}
\end{equation}
where 
\begin{equation}
c_+=\lim_{l\rightarrow\infty} \frac{|A_{l+}|}{|A_{l+1,+}|} ,\:
c_-=\lim_{l\rightarrow\infty} \frac{|A_{l-}|}{|A_{l+1,-}|} ,\:
c_k=\lim_{l\rightarrow\infty} \frac
{\left(A_{lk}^2+\tilde A_{lk}^2\right)^{1/2}}
{\left(A_{l+1,k}^2+\tilde A_{l+1,k}^2\right)^{1/2}}, \;k=1,2.
\label{6ch91}
\end{equation}

\subsection{Analytic functions of a polar 6-compex variable}

The expansion of an analytic function $f(u)$ around $u=u_0$ is
\begin{equation}
f(u)=\sum_{k=0}^\infty \frac{1}{k!} f^{(k)}(u_0)(u-u_0)^k .
\label{6chh91d}
\end{equation}
Since the limit $f^\prime (u_0)=\lim_{u\rightarrow u_0}\{f(u)-f(u_0)\}/(u-u_0)$
is independent of the direction in space along which $u$ is approaching $u_0$,
the function $f(u)$ is said to be analytic, analogously to the case of
functions of regular complex variables. \cite{3} If
$f(u)=\sum_{k=0}^{5}h_kP_k(x_0,x_1,x_2,x_3,x_4,x_{5})$, then
\begin{equation}
\frac{\partial P_0}{\partial x_0} 
=\frac{\partial P_1}{\partial x_1} 
=\frac{\partial P_2}{\partial x_2} 
=\frac{\partial P_3}{\partial x_3}
=\frac{\partial P_4}{\partial x_4}
=\frac{\partial P_5}{\partial x_5}, 
\label{6chh95a}
\end{equation}
\begin{equation}
\frac{\partial P_1}{\partial x_0} 
=\frac{\partial P_2}{\partial x_1} 
=\frac{\partial P_3}{\partial x_2} 
=\frac{\partial P_4}{\partial x_3}
=\frac{\partial P_5}{\partial x_4}
=\frac{\partial P_0}{\partial x_5}, 
\label{6chh95b}
\end{equation}
\begin{equation}
\frac{\partial P_2}{\partial x_0} 
=\frac{\partial P_3}{\partial x_1} 
=\frac{\partial P_4}{\partial x_2} 
=\frac{\partial P_5}{\partial x_3}
=\frac{\partial P_0}{\partial x_4}
=\frac{\partial P_1}{\partial x_5}, 
\label{6chh95c}
\end{equation}
\begin{equation}
\frac{\partial P_3}{\partial x_0} 
=\frac{\partial P_4}{\partial x_1} 
=\frac{\partial P_5}{\partial x_2} 
=\frac{\partial P_0}{\partial x_3}
=\frac{\partial P_1}{\partial x_4}
=\frac{\partial P_2}{\partial x_5}, 
\label{6chh95d}
\end{equation}
\begin{equation}
\frac{\partial P_4}{\partial x_0} 
=\frac{\partial P_5}{\partial x_1} 
=\frac{\partial P_0}{\partial x_2} 
=\frac{\partial P_1}{\partial x_3}
=\frac{\partial P_2}{\partial x_4}
=\frac{\partial P_3}{\partial x_5}, 
\label{6chh95e}
\end{equation}
\begin{equation}
\frac{\partial P_5}{\partial x_0} 
=\frac{\partial P_0}{\partial x_1} 
=\frac{\partial P_1}{\partial x_2} 
=\frac{\partial P_2}{\partial x_3}
=\frac{\partial P_3}{\partial x_4}
=\frac{\partial P_4}{\partial x_5}, 
\label{6chh95f}
\end{equation}
and
\begin{eqnarray}
\lefteqn{\frac{\partial^2 P_k}{\partial x_0\partial x_l}
=\frac{\partial^2 P_k}{\partial x_1\partial x_{l-1}}
=\cdots=
\frac{\partial^2 P_k}{\partial x_{[l/2]}\partial x_{l-[l/2]}}}\nonumber\\
&&=\frac{\partial^2 P_k}{\partial x_{l+1}\partial x_{5}}
=\frac{\partial^2 P_k}{\partial x_{l+2}\partial x_{4}}
=\cdots
=\frac{\partial^2 P_k}{\partial x_{l+1+[(4-l)/2]}
\partial x_{5-[(4-l)/2]}} ,
\label{6ch96}
\end{eqnarray}
for $k,l=0,...,5$.
In Eq. (\ref{6ch96}), $[a]$ denotes the integer part of $a$,
defined as $[a]\leq a<[a]+1$.
In this work, brackets larger than the regular brackets
$[\;]$ do not have the meaning of integer part.

\subsection{Integrals of polar 6-complex functions}

If $f(u)$ is an analytic 6-complex function,
then
\begin{equation}
\oint_\Gamma \frac{f(u)du}{u-u_0}=
2\pi f(u_0)\left[\tilde e_1 
\;{\rm int}(u_{0\xi_1\eta_1},\Gamma_{\xi_1\eta_1})
+\tilde e_2 
\;{\rm int}(u_{0\xi_2\eta_2},\Gamma_{\xi_2\eta_2})
 \right],
\label{6ch120}
\end{equation}
where
\begin{equation}
{\rm int}(M,C)=\left\{
\begin{array}{l}
1 \;\:{\rm if} \;\:M \;\:{\rm is \;\:an \;\:interior \;\:point \;\:of} \;\:C
,\\  
0 \;\:{\rm if} \;\:M \;\:{\rm is \;\:exterior \;\:to}\:\; C ,\\
\end{array}\right.,
\label{6ch118}
\end{equation}
and $u_{0\xi_k\eta_k}$ and $\Gamma_{\xi_k\eta_k}$ are respectively the
projections of the pole $u_0$ and of 
the loop $\Gamma$ on the plane defined by the axes $\xi_k$ and $\eta_k$,
$k=1,2$. 

\subsection{Factorization of 6-complex polynomials}

A polynomial of degree $m$ of the 6-complex variable $u$ has the form
\begin{equation}
P_m(u)=u^m+a_1 u^{m-1}+\cdots+a_{m-1} u +a_m ,
\label{6ch125}
\end{equation}
where $a_l$, for $l=1,...,m$, are 6-complex constants.
If $a_l=\sum_{p=0}^{5}h_p a_{lp}$, and with the notations of Eqs.
(\ref{6ch88a})-(\ref{6ch88c}) applied for $l= 1, \cdots, m$, the polynomial
$P_m(u)$ can be written as
\begin{eqnarray}
\lefteqn{P_m= 
e_+\left(v_+^m +\sum_{l=1}^{m}A_{l+}v_+^{m-l} \right)
+e_-\left(v_-^m +\sum_{l=1}^{m}A_{l-}v_-^{m-l} \right) \nonumber}\\
&&+\sum_{k=1}^{2}
\left[(e_k v_k+\tilde e_k\tilde v_k)^m+
\sum_{l=1}^m(e_k A_{lk}+\tilde e_k\tilde A_{lk})
(e_k v_k+\tilde e_k\tilde v_k)^{m-l} 
\right],
\label{6ch126a}
\end{eqnarray}
where the constants $A_{l+}, A_{l-}, A_{lk}, \tilde A_{lk}$ 
are real numbers.

The polynomial $P_m(u)$ can be written, as 
\begin{eqnarray}
P_m(u)=\prod_{p=1}^m (u-u_p) ,
\label{6ch128c}
\end{eqnarray}
where
\begin{eqnarray}
u_p=e_+ v_{p+}+e_-v_{p-}
+\left(e_1 v_{1p}+\tilde e_1\tilde v_{1p}\right)
+\left(e_2 v_{2p}+\tilde e_2\tilde v_{2p}\right), p=1,...,m.
\label{6ch128d}
\end{eqnarray}
The quantities $v_{p+}$,  $v_{p-}$,  
$e_k v_{kp}+\tilde e_k\tilde v_{kp}$, $p=1,...,m, k=1,2$,
are the roots of the corresponding polynomial in Eq. (\ref{6ch126a}).
The roots $v_{p+}$,  $v_{p-}$ appear in complex-conjugate pairs, and 
$v_{kp}, \tilde v_{kp}$ are real numbers.
Since all these roots may be ordered arbitrarily, the polynomial $P_m(u)$ can
be written in many different ways as a product of linear factors. 

If $P(u)=u^2-1$, the degree is $m=2$, the coefficients of the polynomial are
$a_1=0, a_2=-1$, the coefficients defined in Eqs. (\ref{6ch88a})-(\ref{6ch88c})
are $A_{2+}=-1, A_{2-}=-1, A_{21}=-1, \tilde A_{21}=0,
A_{22}=-1, \tilde A_{22}=0$. The expression of $P(u)$, Eq. (\ref{6ch126a}), is
$v_+^2-e_++v_-^2-e_- +(e_1v_1+\tilde e_1\tilde v_1)^2-e_1+
(e_2v_2+\tilde e_2\tilde v_2)^2-e_2 $. 
The factorization of $P(u)$, Eq. (\ref{6ch128c}), is
$P(u)=(u-u_1)(u-u_2)$, where the roots are
$u_1=\pm e_+\pm e_-\pm e_1\pm  e_2, u_2=-u_1$. If $e_+, e_-, e_1, e_2$ 
are expressed with the aid of Eq. (\ref{6che11}) in terms of $h_1, h_2, h_3,
h_4, h_5$, the factorizations of $P(u)$ are obtained as
\begin{eqnarray}
\lefteqn{\begin{array}{l}
u^2-1=(u+1)(u-1),\\
u^2-1=\left[u+\frac{1}{3}(1+h_1+h_2-2h_3+h_4+h_5)\right]
\left[u-\frac{1}{3}(1+h_1+h_2-2h_3+h_4+h_5)\right],\\
u^2-1=\left[u+\frac{1}{3}(1-h_1+h_2+2h_3+h_4-h_5)\right]
\left[u-\frac{1}{3}(1-h_1+h_2+2h_3+h_4-h_5)\right],\\
u^2-1=\left[u+\frac{1}{3}(2+h_1-h_2+h_3-h_4+h_5)\right]
\left[u-\frac{1}{3}(2+h_1-h_2+h_3-h_4+h_5)\right],\\
u^2-1=\left[u+\frac{1}{3}(-1+2h_2+2h_4)\right]
\left[u-\frac{1}{3}(-1+2h_2+2h_4)\right],\\
u^2-1=\left[u+\frac{1}{3}(2h_1-h_3+2h_5)\right]
\left[u-\frac{1}{3}(2h_1-h_3+2h_5)\right],\\
u^2-1=(u+h_3)(u-h_3),\\
u^2-1=\left[u+\frac{1}{3}(-2+h_1+h_2+h_3+h_4+h_5)\right]
\left[u-\frac{1}{3}(-2+h_1+h_2+h_3+h_4+h_5)\right].
\end{array}\nonumber}\\
&&
\label{ch-factorization}
\end{eqnarray}
It can be checked that 
$(\pm e_+\pm e_-\pm e_1\pm e_2)^2=e_++e_-+e_1+e_2=1$.

\subsection{Representation of polar 6-complex numbers by irreducible matrices}

If the unitary matrix which appears in the expression, 
Eq. (\ref{6ch9ee}), of the variables $\xi_+, \xi_-, 
\xi_1, \eta_1, \xi_k, \eta_k$ in
terms of $x_0, x_1, x_2, x_3, x_4, x_5$ is called $T$,
the irreducible representation of the hypercomplex number $u$ is
\begin{equation}
T U T^{-1}=\left(
\begin{array}{ccccc}
v_+     &     0     &     0   &    0   \\
0       &     v_-   &     0   &    0   \\
0       &     0     &     V_1 &    0   \\
0       &     0     &     0   &    V_2\\
\end{array}
\right),
\label{6ch129a}
\end{equation}
where $U$ is the matrix in Eq. (\ref{6ch24b}),
and $V_k$ are the matrices
\begin{equation}
V_k=\left(
\begin{array}{cc}
v_k           &     \tilde v_k   \\
-\tilde v_k   &     v_k          \\
\end{array}\right),\;\; k=1,2.
\label{6ch130}
\end{equation}

\section{Planar complex numbers in 6 dimensions}

\subsection{Operations with planar complex numbers in 6 dimensions}

The planar hypercomplex number $u$ in 6 dimensions 
is represented as  
\begin{equation}
u=x_0+h_1x_1+h_2x_2+h_3x_3+h_4x_4+h_5x_5. 
\label{6c1a}
\end{equation}
The multiplication rules for the bases 
$h_1, h_2, h_3, h_4, h_5 $ are 
\begin{eqnarray}
\lefteqn{h_1^2=h_2,\;h_2^2=h_4,\;h_3^2=1,\;h_4^2=-h_2,\;h_5^2=-h_4,\;
h_1h_2=h_3,\;h_1h_3=h_4,\;h_1h_4=h_5,\;h_1h_5=-1,\nonumber}\\
&&\;h_2h_3=h_5,\;h_2h_4=-1,\;h_2h_5=-h_1,\;h_3h_4=-h_1,\;h_3h_5=-h_2,
\;h_4h_5=-h_3.
\label{6c1}
\end{eqnarray}
The significance of the composition laws in Eq.
(\ref{6c1}) can be understood by representing the bases 
$1,h_1,h_2,h_3,h_4,h_5$
by points on a circle at the angles $\alpha_k=\pi k/6$.
The product $h_j h_k$ will be represented by the point of the circle at the
angle $\pi (j+k)/12$, $j,k=0,1,...,5$. If $\pi\leq\pi (j+k)/12\leq 2\pi$, the
point is opposite to the basis $h_l$ of angle $\alpha_l=\pi (j+k)/6-\pi$. \\

The sum of the 6-complex numbers $u$ and $u^\prime$ is
\begin{equation}
u+u^\prime=x_0+x^\prime_0+h_1(x_1+x^\prime_1)+h_1(x_2+x^\prime_2)
+h_3(x_3+x^\prime_3)+h_4(x_4+x^\prime_4)+h_5(x_5+x^\prime_5).
\label{6c2}
\end{equation}
The product of the numbers $u, u^\prime$ is
\begin{equation}
\begin{array}{l}
uu^\prime=x_0 x_0^\prime -x_1x_5^\prime-x_2 x_4^\prime-x_3x_3^\prime
-x_4x_2^\prime-x_5 x_1^\prime\\
+h_1(x_0 x_1^\prime+x_1x_0^\prime-x_2x_5^\prime-x_3x_4^\prime
-x_4 x_3^\prime-x_5 x_2^\prime) \\
+h_2(x_0 x_2^\prime+x_1x_1^\prime+x_2x_0^\prime-x_3x_5^\prime
-x_4 x_4^\prime-x_5 x_3^\prime) \\
+h_3(x_0 x_3^\prime+x_1x_2^\prime+x_2x_1^\prime+x_3x_0^\prime
-x_4 x_5^\prime-x_5 x_4^\prime) \\
+h_4(x_0 x_4^\prime+x_1x_3^\prime+x_2x_2^\prime+x_3x_1^\prime
+x_4 x_0^\prime-x_5 x_5^\prime) \\
+h_5(x_0 x_5^\prime+x_1x_4^\prime+x_2x_3^\prime+x_3x_2^\prime
+x_4 x_1^\prime+x_5 x_0^\prime).
\end{array}
\label{6c3}
\end{equation}

The relation between the variables $v_1,\tilde v_1, v_2, \tilde v_2,
v_3,\tilde v_3$ and $x_0,x_1,x_2,x_3,x_4,x_5$ are
\begin{equation}
\left(
\begin{array}{c}
v_1\\
\tilde v_1\\
v_2\\
\tilde v_2\\
v_3\\
\tilde v_3\\
\end{array}\right)
=\left(
\begin{array}{cccccc}
1&\frac{\sqrt{3}}{2}&\frac{1}{2}&0&-\frac{1}{2}&-\frac{\sqrt{3}}{2}\\
0&\frac{1}{2}&\frac{\sqrt{3}}{2}&1&\frac{\sqrt{3}}{2}&\frac{1}{2}\\
1&0 &-1 &0 &1 &0 \\
0&1&0 &-1 &0 &1 \\
1&-\frac{\sqrt{3}}{2}&\frac{1}{2}&0&-\frac{1}{2}&\frac{\sqrt{3}}{2}\\
0&\frac{1}{2}&-\frac{\sqrt{3}}{2}&1&-\frac{\sqrt{3}}{2}&\frac{1}{2}\\
\end{array}
\right)
\left(
\begin{array}{c}
x_0\\
x_1\\
x_2\\
x_3\\
x_4\\
x_5
\end{array}
\right).
\label{6c9e}
\end{equation}
The other variables are $v_4=v_3, \tilde v_4=-\tilde v_3,
v_5=v_2, \tilde v_5=-\tilde v_2, v_6=v_1, \tilde v_6=-\tilde v_1$. 
The variables $v_1, \tilde v_1, v_2, \tilde v_2, v_3, \tilde v_3$ will be
called canonical planar 6-complex variables.

\subsection{Geometric representation of planar complex numbers in 6 dimensions}

The 6-complex number $u=x_0+h_1x_1+h_2x_2+h_3x_3+h_4x_4+h_5x_5$
is represented by 
the point $A$ of coordinates $(x_0,x_1,x_2,x_3,x_4,x_5)$. 
The distance from the origin $O$ of the 6-dimensional space to the point $A$
has the expression 
\begin{equation}
d^2=x_0^2+x_1^2+x_2^2+x_3^2+x_4^2+x_5^2,
\label{6c10}
\end{equation}
is called modulus of the 6-complex number 
$u$, and is designated by $d=|u|$.
The modulus has the property that
\begin{equation}
|u^\prime u^{\prime\prime}|\leq \sqrt{3}|u^\prime||u^{\prime\prime}| .
\label{6c79}
\end{equation}

The exponential and trigonometric forms of the 6-complex number $u$ can be
obtained conveniently in a rotated system of axes defined by a transformation
which has the form
\begin{equation}
\left(
\begin{array}{c}
\xi_1\\
\tilde \xi_1\\
\xi_2\\
\tilde \xi_2\\
\xi_3\\
\tilde \xi_3\\
\end{array}\right)
=\left(
\begin{array}{cccccc}
\frac{1}{\sqrt{3}}&\frac{1}{2}&\frac{1}{2\sqrt{3}}&0&
-\frac{1}{2\sqrt{3}}&-\frac{1}{2}\\
0&\frac{1}{2\sqrt{3}}&\frac{1}{2}&\frac{1}{\sqrt{3}}&
\frac{1}{2}&\frac{1}{2\sqrt{3}}\\
\frac{1}{\sqrt{3}}&0 &-\frac{1}{\sqrt{3}} &0 &\frac{1}{\sqrt{3}} &0 \\
0&\frac{1}{\sqrt{3}}&0 &-\frac{1}{\sqrt{3}} &0 &\frac{1}{\sqrt{3}} \\
\frac{1}{\sqrt{3}}&-\frac{1}{2}&\frac{1}{2\sqrt{3}}&0&
-\frac{1}{2\sqrt{3}}&\frac{1}{2}\\
0&\frac{1}{2\sqrt{3}}&-\frac{1}{2}&\frac{1}{\sqrt{3}}&
-\frac{1}{2}&\frac{1}{2\sqrt{3}}\\
\end{array}
\right)
\left(
\begin{array}{c}
x_0\\
x_1\\
x_2\\
x_3\\
x_4\\
x_5
\end{array}
\right).
\label{6c9ee}
\end{equation}

The lines of the matrices in Eq. (\ref{6c9ee}) give the components
of the 6 vectors of the new basis system of axes. These vectors have unit
length and are orthogonal to each other.
The relations between the two sets of variables are
\begin{equation}
v_k= \sqrt{3}\xi_k , \tilde v_k= \sqrt{3}\eta_k, 
\label{6c12b}
\end{equation}
for $k=1,2,3$.
 
The radius $\rho_k$ and the azimuthal angle $\phi_k$ in the plane of the axes
$v_k,\tilde v_k$ are
\begin{equation}
\rho_k^2=v_k^2+\tilde v_k^2, \:\cos\phi_k=v_k/\rho_k,
\:\sin\phi_k=\tilde v_k/\rho_k, 
\label{6c19a}
\end{equation}
where $0\leq \phi_k<2\pi ,  \;k=1,2,3$,
so that there are 3 azimuthal angles.
The planar angles $\psi_{k-1}$ are
\begin{equation}
\tan\psi_1=\rho_1/\rho_2,  \;\tan\psi_2=\rho_1/\rho_3, 
\label{6c19b}
\end{equation}
where $0\leq\psi_1\leq\pi/2,\;0\leq\psi_2\leq\pi/2$, 
so that there are 2 planar angles.
The amplitude of an 6-complex number $u$ is
\begin{equation}
\rho=\left(\rho_1\rho_2\rho_3\right)^{1/3}.
\label{6c50aa}
\end{equation}
It can be checked that
\begin{equation}
d^2=\frac{1}{3}(\rho_1^2+\rho_2^2+\rho_3^2).
\label{6c13}
\end{equation}

If $u=u^\prime u^{\prime\prime}$, the parameters of the hypercomplex numbers
are related by
\begin{equation}
\rho_k=\rho_k^\prime\rho_k^{\prime\prime}, 
\label{6c21b}
\end{equation}
\begin{equation}
\tan\psi_k=\tan\psi_k^\prime \tan\psi_k^{\prime\prime},  
\label{6c21d}
\end{equation}
\begin{equation}
\phi_k=\phi_k^\prime+\phi_k^{\prime\prime},
\label{6c21e}
\end{equation}
\begin{equation}
v_k=v_k^\prime v_k^{\prime\prime}-\tilde v_k^\prime 
\tilde v_k^{\prime\prime},\;
\tilde v_k=v_k^\prime \tilde v_k^{\prime\prime}
+\tilde v_k^\prime v_k^{\prime\prime},
\label{6c22}
\end{equation}
\begin{equation}
\rho=\rho^\prime\rho^{\prime\prime} ,
\label{6c24}
\end{equation}
where $k=1,2,3$.

The 6-complex planar number
$u=x_0+h_1x_1+h_2x_2+h_3x_3+h_4x_4+h_5x_5$ can be  represented by the matrix
\begin{equation}
U=\left(
\begin{array}{cccccc}
x_0      &    x_1     &    x_2   &    x_3  &  x_4   & x_5\\
-x_5     &    x_0     &    x_1   &    x_2  &  x_3   & x_4\\
-x_4     &   -x_5     &    x_0   &    x_1  &  x_2   & x_3\\
-x_3     &   -x_4     &   -x_5   &    x_0  &  x_1   & x_2\\
-x_2     &   -x_3     &   -x_4   &   -x_5  &  x_0   & x_1\\
-x_1     &   -x_2     &   -x_3   &   -x_4  &  -x_5  & x_0\\
\end{array}
\right).
\label{6c24b}
\end{equation}
The product $u=u^\prime u^{\prime\prime}$ is
represented by the matrix multiplication $U=U^\prime U^{\prime\prime}$.

\subsection{The planar 6-dimensional cosexponential functions}

The planar cosexponential functions in 6 dimensions are
\begin{equation}
f_{6k}(y)=\sum_{p=0}^\infty (-1)^p \frac{y^{k+6p}}{(k+6p)!}, 
\label{6c29}
\end{equation}
for $k=0,...,5$.
The planar cosexponential functions  of even index $k$ are
even functions, $f_{6,2l}(-y)=f_{6,2l}(y)$,  
and the planar cosexponential functions of odd index 
are odd functions, $f_{6,2l+1}(-y)=-f_{6,2l+1}(y)$, $l=0,1,2$. 
The exponential function of the quantity $h_k y$ is
\begin{equation}
\begin{array}{l}
e^{h_1 y}=f_{60}(y)+h_1f_{61}(y)+h_2f_{62}(y)
+h_3f_{63}(y)+h_4f_{64}(y)+h_5f_{65}(y),\\
e^{h_2 y}=g_{60}(y)-g_{63}(y)+h_2\{g_{61}(y)
-g_{64}(y)\}+h_4\{g_{62}(y)-g_{65}(y)\},\\
e^{h_3 y}=f_{60}(y)-f_{62}(y)+f_{64}(y)+h_3\{f_{61}(y)-f_{63}(y)+f_{65}(y)\},\\
e^{h_4 y}=g_{60}(y)+g_{63}(y)-h_2\{g_{62}(y)
+g_{65}(y)\}+h_4\{g_{61}(y)+g_{64}(y)\},\\
e^{h_5 y}=f_{60}(y)+h_1f_{65}(y)-h_2f_{64}(y)
+h_3f_{63}(y)-h_4f_{62}(y)+h_5f_{61}(y).\\
\end{array}
\label{6c28b}
\end{equation}
The relations for $h_2$ and $h_4$ can be written equivalently as 
$e^{h_2 y}=f_{30}+h_2f_{31}+h_4 f_{32}, e^{h_4 y}=g_{30}-h_2f_{32}+h_4 g_{31}$,
and the relation for $h_3$ can be written as $e^{h_3 y}=f_{20}+h_3f_{21}$,
which is the same as $e^{h_3 y}=\cos y+h_3\sin y$. 

The planar 6-dimensional cosexponential functions $f_{6k}(y)$ are related to
the polar 6-dimensional cosexponential function $g_{6k}(y)$ by the relations
\begin{equation}
f_{6k}(y)=e^{-i\pi k/6}g_{6k}\left(e^{i\pi/6}y\right), 
\label{6c30a}
\end{equation}
for $k=0,...,5$.
The planar 6-dimensional cosexponential functions $f_{6k}(y)$ are related to
the polar 6-dimensional cosexponential function $g_{6k}(y)$ also
by the relations 
\begin{equation}
f_{6k}(y)=e^{-i\pi k/2}g_{6k}(iy), 
\label{6c30ax}
\end{equation}
for $k=0,...,5$.
The expressions of the planar 6-dimensional cosexponential functions are
\begin{equation}
\begin{array}{l}
f_{60}(y)=\frac{1}{3}\cos y
+\frac{2}{3}\cosh\frac{\sqrt{3}}{2}y\cos\frac{y}{2},\\
f_{61}(y)=\frac{1}{3}\sin y
+\frac{\sqrt{3}}{3}\sinh\frac{\sqrt{3}}{2}y\cos\frac{y}{2}
+\frac{1}{3}\cosh\frac{\sqrt{3}}{2}y\sin\frac{y}{2},\\
f_{62}(y)=-\frac{1}{3}\cos y+\frac{1}{3}\cosh\frac{\sqrt{3}}{2}y\cos\frac{y}{2}
+\frac{\sqrt{3}}{3}\sinh\frac{\sqrt{3}}{2}y\sin\frac{y}{2},\\
f_{63}(y)=-\frac{1}{3}\sin y
+\frac{2}{3}\cosh\frac{\sqrt{3}}{2}y\sin\frac{y}{2},\\
f_{64}(y)=\frac{1}{3}\cos y-\frac{1}{3}\cosh\frac{\sqrt{3}}{2}y\cos\frac{y}{2}
+\frac{\sqrt{3}}{3}\sinh\frac{\sqrt{3}}{2}y\sin\frac{y}{2},\\
f_{65}(y)=\frac{1}{3}\sin y
-\frac{\sqrt{3}}{3}\sinh\frac{\sqrt{3}}{2}y\cos\frac{y}{2}
+\frac{1}{3}\cosh\frac{\sqrt{3}}{2}y\sin\frac{y}{2}.\\
\end{array}
\label{6c30x}
\end{equation}
The planar 6-dimensional cosexponential functions can be written as
\begin{equation}
f_{6k}(y)=\frac{1}{6}\sum_{l=1}^{6}
\exp\left[y\cos\left(\frac{\pi (2l-1)}{6}\right)
\right]
\cos\left[y\sin\left(\frac{\pi (2l-1)}{6}\right)-\frac{\pi (2l-1)k}{6}\right], 
\label{6c30}
\end{equation}
for $k=0,...,5$.
The graphs of the planar 6-dimensional cosexponential functions are shown in
Fig. 4. 

It can be checked that
\begin{equation}
\sum_{k=0}^{5}f_{6k}^2(y)=\frac{1}{3}+\frac{2}{3}\cosh\sqrt{3}y.
\label{6c34a}
\end{equation}

The addition theorems for the planar 6-dimensional cosexponential functions are
\begin{eqnarray}
\lefteqn{\begin{array}{l}
g_{60}(y+z)=g_{60}(y)g_{60}(z)-g_{61}(y)g_{65}(z)-g_{62}(y)g_{64}(z)
-g_{63}(y)g_{63}(z)-g_{64}(y)g_{62}(z)-g_{65}(y)g_{61}(z) ,\\
g_{61}(y+z)=g_{60}(y)g_{61}(z)+g_{61}(y)g_{60}(z)-g_{62}(y)g_{65}(z)
-g_{63}(y)g_{64}(z)-g_{64}(y)g_{63}(z)-g_{65}(y)g_{62}(z) ,\\
g_{62}(y+z)=g_{60}(y)g_{62}(z)+g_{61}(y)g_{61}(z)+g_{62}(y)g_{60}(z)
-g_{63}(y)g_{65}(z)-g_{64}(y)g_{64}(z)-g_{65}(y)g_{63}(z) ,\\
g_{63}(y+z)=g_{60}(y)g_{63}(z)+g_{61}(y)g_{62}(z)+g_{62}(y)g_{61}(z)
+g_{63}(y)g_{60}(z)-g_{64}(y)g_{65}(z)-g_{65}(y)g_{64}(z) ,\\
g_{64}(y+z)=g_{60}(y)g_{64}(z)+g_{61}(y)g_{63}(z)+g_{62}(y)g_{62}(z)
+g_{63}(y)g_{61}(z)+g_{64}(y)g_{60}(z)-g_{65}(y)g_{65}(z) ,\\
g_{65}(y+z)=g_{60}(y)g_{65}(z)+g_{61}(y)g_{64}(z)+g_{62}(y)g_{63}(z)
+g_{63}(y)g_{62}(z)+g_{64}(y)g_{61}(z)+g_{65}(y)g_{60}(z) .
\end{array}\nonumber}\\
&&
\label{6c35a}
\end{eqnarray}
It can be shown that
\begin{equation}
\begin{array}{l}
\{f_{60}(y)+h_1f_{61}(y)+h_2f_{62}(y)+h_3f_{63}(y)
+h_4f_{64}(y)+h_5f_{65}(y)\}^l\\
\hspace*{0.5cm}=f_{60}(ly)+h_1f_{61}(ly)
+h_2f_{62}(ly)+h_3f_{63}(ly)+h_4f_{64}(ly)+h_5f_{65}(ly),\\
\{g_{60}(y)-g_{63}(y)+h_2\{g_{61}(y)
-g_{64}(y)\}+h_4\{g_{62}(y)-g_{65}(y)\}\}^l\\
\hspace*{0.5cm}=g_{60}(ly)-g_{63}(ly)
+h_2\{g_{61}(ly)-g_{64}(ly)\}+h_4\{g_{62}(ly)-g_{65}(ly)\},\\
\{f_{60}(y)-f_{62}(y)+f_{64}(y)+h_3\{f_{61}(y)-f_{63}(y)+f_{65}(y)\}\}^l\\
\hspace*{0.5cm}=f_{60}(ly)-f_{62}(ly)
+f_{64}(ly)+h_3\{f_{61}(ly)-f_{63}(ly)+f_{65}(ly)\},\\
\{g_{60}(y)+g_{63}(y)-h_2\{g_{62}(y)
+g_{65}(y)\}+h_4\{g_{61}(y)+g_{64}(y)\}\}^l\\
\hspace*{0.5cm}=g_{60}(ly)+g_{63}(ly)
-h_2\{g_{62}(ly)+g_{65}(ly)\}+h_4\{g_{61}(ly)+g_{64}(ly)\},\\
\{f_{60}(y)+h_1f_{65}(y)-h_2f_{64}(y)
+h_3f_{63}(y)-h_4f_{62}(y)+h_5f_{61}(y)\}^l\\
\hspace*{0.5cm}=f_{60}(ly)+h_1f_{65}(ly)
-h_2f_{64}(ly)+h_3f_{63}(ly)-h_4f_{62}(ly)+h_5f_{61}(ly).\\
\end{array}
\label{6c37b}
\end{equation}

The derivatives of the planar cosexponential functions
are related by
\begin{equation}
\frac{df_{60}}{du}=-f_{65}, \:
\frac{df_{61}}{du}=f_{60}, \:
\frac{df_{62}}{du}=f_{61}, \:
\frac{df_{63}}{du}=f_{62}, \:
\frac{df_{64}}{du}=f_{63}, \:
\frac{df_{65}}{du}=f_{64}.
\label{6c45}
\end{equation}

\subsection{Exponential and trigonometric forms of planar 6-complex numbers}

The exponential and trigonometric forms of planar 6-complex
numbers can be expressed with the aid of the hypercomplex bases 
\begin{equation}
\left(
\begin{array}{c}
e_1\\
\tilde e_1\\
e_2\\
\tilde e_2\\
e_3\\
\tilde e_3\\
\end{array}\right)
=\left(
\begin{array}{cccccc}
\frac{1}{3}&\frac{\sqrt{3}}{6}&\frac{1}{6}&0&-\frac{1}{6}&-\frac{\sqrt{3}}{6}\\
0&\frac{1}{6}&\frac{\sqrt{3}}{6}&\frac{1}{3}&\frac{\sqrt{3}}{6}&\frac{1}{6}\\
\frac{1}{3}&0 &-\frac{1}{3} &0 &\frac{1}{3} &0 \\
0&\frac{1}{3}&0 &-\frac{1}{3} &0 &\frac{1}{3} \\
\frac{1}{3}&-\frac{\sqrt{3}}{6}&\frac{1}{6}&0&-\frac{1}{6}&\frac{\sqrt{3}}{6}\\
0&\frac{1}{6}&-\frac{\sqrt{3}}{6}&\frac{1}{3}&-\frac{\sqrt{3}}{6}&\frac{1}{6}\\
\end{array}
\right)
\left(
\begin{array}{c}
1\\
h_1\\
h_2\\
h_3\\
h_4\\
h_5
\end{array}
\right).
\label{6ce11}
\end{equation}

The multiplication relations for the bases $e_k, \tilde e_k$ are
\begin{eqnarray}
e_k^2=e_k, \tilde e_k^2=-e_k, e_k \tilde e_k=\tilde e_k , e_ke_l=0, e_k\tilde
e_l=0, \tilde e_k\tilde e_l=0, \;k,l=1,2,3, \;k\not=l.
\label{6ce12a}
\end{eqnarray}
The moduli of the bases $e_k, \tilde e_k$ are
\begin{equation}
|e_k|=\sqrt{\frac{1}{3}}, |\tilde e_k|=\sqrt{\frac{1}{3}}, 
\label{6ce12c}
\end{equation}
for $k=1,2,3$.
It can be shown that
\begin{eqnarray}
x_0+h_1x_1+h_2x_2+h_3x_3+h_4x_4+h_5x_5
= \sum_{k=1}^3 (e_k v_k+\tilde e_k \tilde v_k).
\label{6ce13a}
\end{eqnarray}
The ensemble $e_1, \tilde e_1, e_2, 
\tilde e_2, e_3, \tilde e_3$ will be called the
canonical planar 6-complex base, and Eq. (\ref{6ce13a}) gives the canonical
form of the planar 6-complex number.

The exponential form of the 6-complex number $u$ is
\begin{eqnarray}
\lefteqn{u=\rho\exp\left\{\frac{1}{3}(h_2-h_4)\ln\tan\psi_1
+\frac{1}{6}(\sqrt{3}h_1-h_2+h_4-\sqrt{3}h_5)\ln\tan\psi_2\right.\nonumber}\\
&&\left.+\tilde e_1\phi_1+\tilde e_2\phi_2+\tilde e_3\phi_3\right\}.
\label{6c50a}
\end{eqnarray}
The
trigonometric form of the 6-complex number $u$ is
\begin{eqnarray}
\lefteqn{u=d
\sqrt{3}
\left(1+\frac{1}{\tan^2\psi_1}+\frac{1}{\tan^2\psi_2}\right)^{-1/2}\nonumber}\\
&&\left(e_1+\frac{e_2}{\tan\psi_1}+\frac{e_3}{\tan\psi_2}\right)
\exp\left(\tilde e_1\phi_1+\tilde e_2\phi_2+\tilde e_3\phi_3\right).
\label{6c52a}
\end{eqnarray}
The modulus $d$ and the amplitude $\rho$ are related by
\begin{eqnarray}
d=\rho \frac{2^{1/3}}{\sqrt6}
\left(\tan\psi_1\tan\psi_2\right)^{1/3}
\left(1+\frac{1}{\tan^2\psi_1}+\frac{1}{\tan^2\psi_2}\right)^{1/2}.
\label{6c53a}
\end{eqnarray}

\subsection{Elementary functions of a planar 6-complex variable}

The logarithm and power functions of the 6-complex number $u$ exist for all
$x_0,...,x_5$ and are
\begin{eqnarray}
\lefteqn{\ln u=\ln \rho+
\frac{1}{3}(h_2-h_4)\ln\tan\psi_1
+\frac{1}{6}(\sqrt{3}h_1-h_2+h_4-\sqrt{3}h_5)\ln\tan\psi_2\nonumber}\\
&&+\tilde e_1\phi_1+\tilde e_2\phi_2+\tilde e_3\phi_3 ,
\label{6c56a}
\end{eqnarray}
\begin{equation}
u^m=\sum_{k=1}^{3}
\rho_k^m(e_k\cos m\phi_k+\tilde e_k\sin m\phi_k).
\label{6c59a}
\end{equation}

The exponential of the 6-complex variable $u$ is
\begin{eqnarray}
e^u= 
\sum_{k=1}^{3}e^{v_k}\left(e_k \cos \tilde v_k+\tilde e_k \sin\tilde
v_k\right).
\label{6c73a}
\end{eqnarray}

The trigonometric functions of the
6-complex variable $u$ are
\begin{equation}
\cos u=\sum_{k=1}^{3}\left(e_k \cos v_k\cosh \tilde v_k
-\tilde e_k \sin v_k\sinh\tilde v_k\right),
\label{6c74a}
\end{equation}
\begin{equation}
\sin u= 
\sum_{k=1}^{3}\left(e_k \sin v_k\cosh \tilde v_k
+\tilde e_k \cos v_k\sinh\tilde v_k\right).
\label{6c74b}
\end{equation}
The hyperbolic functions of the
6-complex variable $u$ are
\begin{equation}
\cosh u=
\sum_{k=1}^{3}\left(e_k \cosh v_k\cos \tilde v_k
+\tilde e_k \sinh v_k\sin\tilde v_k\right),
\label{6c75a}
\end{equation}
\begin{equation}
\sinh u=
\sum_{k=1}^{3}\left(e_k \sinh v_k\cos \tilde v_k
+\tilde e_k \cosh v_k\sin\tilde v_k\right).
\label{6c75b}
\end{equation}

\subsection{Power series of 6-complex numbers}

A power series of the 6-complex variable $u$ is a series of the form
\begin{equation}
a_0+a_1 u + a_2 u^2+\cdots +a_l u^l+\cdots .
\label{6c83}
\end{equation}
Since
\begin{equation}
|au^l|\leq 3^{l/2} |a| |u|^l ,
\label{6c82}
\end{equation}
the series is absolutely convergent for 
\begin{equation}
|u|<c,
\label{6c86}
\end{equation}
where 
\begin{equation}
c=\lim_{l\rightarrow\infty} \frac{|a_l|}{\sqrt{3}|a_{l+1}|} .
\label{6c87}
\end{equation}

If $a_l=\sum_{p=0}^{5} h_p a_{lp}$, and
\begin{equation}
A_{lk}=\sum_{p=0}^{5} a_{lp}\cos\frac{\pi (2k-1)p}{6},
\label{6c88b}
\end{equation}
\begin{equation}
\tilde A_{lk}=\sum_{p=0}^{5} a_{lp}\sin\frac{\pi (2k-1)p}{6},
\label{6c88c}
\end{equation}
where $k=1,2,3$, the series (\ref{6c83}) can be written as
\begin{equation}
\sum_{l=0}^\infty \left[
\sum_{k=1}^{3}
(e_k A_{lk}+\tilde e_k\tilde A_{lk})(e_k v_k+\tilde e_k\tilde v_k)^l 
\right].
\label{6c89a}
\end{equation}
The series is absolutely convergent for   
\begin{equation}
\rho_k<c_k, k=1,2,3,
\label{6c90}
\end{equation}
where 
\begin{equation}
c_k=\lim_{l\rightarrow\infty} \frac
{\left[A_{lk}^2+\tilde A_{lk}^2\right]^{1/2}}
{\left[A_{l+1,k}^2+\tilde A_{l+1,k}^2\right]^{1/2}} .
\label{6c91}
\end{equation}

\subsection{Analytic functions of a planar 6-complex variable}

The expansion of an analytic function $f(u)$ around $u=u_0$ is
\begin{equation}
f(u)=\sum_{k=0}^\infty \frac{1}{k!} f^{(k)}(u_0)(u-u_0)^k .
\label{6ch91d}
\end{equation}

If 
$f(u)=\sum_{k=0}^{5}h_kP_k(x_0,...,x_5)$,
then
\begin{equation}
\frac{\partial P_0}{\partial x_0} 
=\frac{\partial P_1}{\partial x_1} 
=\frac{\partial P_2}{\partial x_2} 
=\frac{\partial P_3}{\partial x_3}
=\frac{\partial P_4}{\partial x_4}
=\frac{\partial P_5}{\partial x_5}, 
\label{6ch95a}
\end{equation}
\begin{equation}
\frac{\partial P_1}{\partial x_0} 
=\frac{\partial P_2}{\partial x_1} 
=\frac{\partial P_3}{\partial x_2} 
=\frac{\partial P_4}{\partial x_3}
=\frac{\partial P_5}{\partial x_4}
=-\frac{\partial P_0}{\partial x_5}, 
\label{6ch95b}
\end{equation}
\begin{equation}
\frac{\partial P_2}{\partial x_0} 
=\frac{\partial P_3}{\partial x_1} 
=\frac{\partial P_4}{\partial x_2} 
=\frac{\partial P_5}{\partial x_3}
=-\frac{\partial P_0}{\partial x_4}
=-\frac{\partial P_1}{\partial x_5}, 
\label{6ch95c}
\end{equation}
\begin{equation}
\frac{\partial P_3}{\partial x_0} 
=\frac{\partial P_4}{\partial x_1} 
=\frac{\partial P_5}{\partial x_2} 
=-\frac{\partial P_0}{\partial x_3}
=-\frac{\partial P_1}{\partial x_4}
=-\frac{\partial P_2}{\partial x_5}, 
\label{6ch95d}
\end{equation}
\begin{equation}
\frac{\partial P_4}{\partial x_0} 
=\frac{\partial P_5}{\partial x_1} 
=-\frac{\partial P_0}{\partial x_2} 
=-\frac{\partial P_1}{\partial x_3}
=-\frac{\partial P_2}{\partial x_4}
=-\frac{\partial P_3}{\partial x_5}, 
\label{6ch95e}
\end{equation}
\begin{equation}
\frac{\partial P_5}{\partial x_0} 
=-\frac{\partial P_0}{\partial x_1} 
=-\frac{\partial P_1}{\partial x_2} 
=-\frac{\partial P_2}{\partial x_3}
=-\frac{\partial P_3}{\partial x_4}
=-\frac{\partial P_4}{\partial x_5}, 
\label{6ch95f}
\end{equation}
and
\begin{eqnarray}
\lefteqn{\frac{\partial^2 P_k}{\partial x_0\partial x_l}
=\frac{\partial^2 P_k}{\partial x_1\partial x_{l-1}}
=\cdots=
\frac{\partial^2 P_k}{\partial x_{[l/2]}\partial x_{l-[l/2]}}}\nonumber\\
&&=-\frac{\partial^2 P_k}{\partial x_{l+1}\partial x_5}
=-\frac{\partial^2 P_k}{\partial x_{l+2}\partial x_4}
=\cdots
=-\frac{\partial^2 P_k}{\partial x_{l+1+[(4-l)/2]}
\partial x_{5-[(4-l)/2]}} .
\label{6c96}
\end{eqnarray}

\subsection{Integrals of planar 6-complex functions}

If $f(u)$ is an analytic 6-complex function,
then
\begin{equation}
\oint_\Gamma \frac{f(u)du}{u-u_0}=
2\pi f(u_0)\left\{
\tilde e_1 \;{\rm int}(u_{0\xi_1\eta_1},\Gamma_{\xi_1\eta_1})
+\tilde e_2 \;{\rm int}(u_{0\xi_2\eta_2},\Gamma_{\xi_2\eta_2})
+\tilde e_3 \;{\rm int}(u_{0\xi_3\eta_3},\Gamma_{\xi_3\eta_3})\right\},
\label{6c120}
\end{equation}
where $u_{0\xi_k\eta_k}$ and $\Gamma_{\xi_k\eta_k}$ are respectively the
projections of the point $u_0$ and of 
the loop $\Gamma$ on the plane defined by the axes $\xi_k$ and $\eta_k$,
$k=1,2,3$. 

\subsection{Factorization of 6-complex polynomials}

A polynomial of degree $m$ of the 6-complex variable $u$ has the form
\begin{equation}
P_m(u)=u^m+a_1 u^{m-1}+\cdots+a_{m-1} u +a_m ,
\label{6c125}
\end{equation}
where $a_l$, for $l=1,...,m$, are 6-complex constants.
If $a_l=\sum_{p=0}^{5}h_p a_{lp}$, and with the
notations of Eqs. (\ref{6c88b})-(\ref{6c88c}) 
applied for $l= 1, \cdots, m$, the
polynomial $P_m(u)$ can be written as 
\begin{eqnarray}
P_m= 
\sum_{k=1}^{3}
\left[(e_k v_k+\tilde e_k\tilde v_k)^m+
\sum_{l=1}^m(e_k A_{lk}+\tilde e_k\tilde A_{lk})
(e_k v_k+\tilde e_k\tilde v_k)^{m-l} 
\right],
\label{6c126a}
\end{eqnarray}
where the constants $A_{lk}, \tilde A_{lk}$ are real numbers.

The polynomial $P_m(u)$ can be written as a product of factors 
\begin{eqnarray}
P_m(u)=\prod_{p=1}^m (u-u_p) ,
\label{6c128c}
\end{eqnarray}
where
\begin{eqnarray}
u_p=\sum_{k=1}^{3}\left(e_k v_{kp}+\tilde e_k\tilde v_{kp}\right), 
\label{6c128d}
\end{eqnarray}
for $p=1,...,m$.
The quantities  
$e_k v_{kp}+\tilde e_k\tilde v_{kp}$, $p=1,...,m, k=1,2,3$,
are the roots of the corresponding polynomial in Eq. (\ref{6c126a}) and are
real numbers.
Since these roots may be ordered arbitrarily, the polynomial $P_m(u)$ can be
written in many different ways as a product of linear factors. 

If $P(u)=u^2+1$, the degree is $m=2$, the coefficients of the polynomial are
$a_1=0, a_2=1$, the coefficients defined in Eqs. (\ref{6c88b})-(\ref{6c88c})
are $A_{21}=1, \tilde A_{21}=0,
A_{22}=1, \tilde A_{22}=0, A_{23}=1, \tilde A_{23}=0
$. The expression, Eq. (\ref{6c126a}), is
P(u)=$(e_1v_1+\tilde e_1\tilde v_1)^2+e_1+
(e_2v_2+\tilde e_2\tilde v_2)^2+e_2+(e_3v_3+\tilde e_3\tilde v_3)^2+e_3 $. 
The factorization of $P(u)$, Eq. (\ref{6c128c}), is
$P(u)=(u-u_1)(u-u_2)$, where the roots are
$u_1=\pm \tilde e_1 \pm \tilde e_2\pm \tilde e_3, u_2=-u_1$. 
If $\tilde e_1, \tilde e_2, \tilde e_3$ 
are expressed with the aid of Eq. (\ref{6ce11}) in terms of $h_1, h_2, h_3,
h_4, h_5$, the factorizations of $P(u)$ are obtained as
\begin{eqnarray}
\lefteqn{\begin{array}{l}
u^2+1=\left[u+\frac{1}{3}(2h_1+h_3+2h_5)\right]
\left[u-\frac{1}{3}(2h_1+h_3+2h_5)\right],\\
u^2+1=\left[u+\frac{1}{3}(h_1+\sqrt{3}h_2-h_3+\sqrt{3}h_4+h_5)\right]
\left[u-\frac{1}{3}(h_1+\sqrt{3}h_2-h_3+\sqrt{3}h_4+h_5)\right],\\
u^2+1=(u+h_3)(u-h_3),\\
u^2+1=\left[u+\frac{1}{3}(-h_1+\sqrt{3}h_2+h_3+\sqrt{3}h_4-h_5)\right]
\left[u-\frac{1}{3}(-h_1+\sqrt{3}h_2+h_3+\sqrt{3}h_4-h_5)\right].
\end{array}\nonumber}\\
&&
\label{c-factorization}
\end{eqnarray}
It can be checked that 
$(\pm \tilde e_1\pm \tilde e_2+\pm \tilde e_3)^2=-e_1-e_2-e_3=-1$.

\subsection{Representation of planar 6-complex numbers by irreducible matrices}

If the unitary matrix written in Eq. (\ref{6c9ee}) is called $T$,
the matric $T U T^{-1}$ provides an irreducible representation 
\cite{4} of the planar
hypercomplex number $u$, 
\begin{equation}
T U T^{-1}=\left(
\begin{array}{ccc}
V_1      &     0   &    0   \\
0        &     V_2 &    0   \\
0        &     0   &    V_3 \\
\end{array}
\right),
\label{6c129}
\end{equation}
where $U$ is the matrix in Eq. (\ref{6c24b}) used to represent the 6-complex
number $u$, and the matrices $V_k$ are
\begin{equation}
V_k=\left(
\begin{array}{cc}
v_k           &     \tilde v_k   \\
-\tilde v_k   &     v_k          \\
\end{array}\right),
\label{6c130}
\end{equation}
for $ k=1,2,3$.

\section{Conclusions}

The operations of addition and multiplication of the polar 6-complex numbers
introduced in this 
work have a geometric interpretation based on the amplitude $\rho$,
the modulus $d$ and the polar, planar and azimuthal angles $\theta_+, \theta_-,
\psi_1, \phi_1, \phi_2$. 
If $v_+>0$ and $v_->0$,
the polar 6-complex numbers can be written in exponential and
trigonometric forms with the aid of the modulus, amplitude and the angular
variables. 
The polar 6-complex functions defined by series of powers are analytic, and 
the partial derivatives of the components of the 6-complex functions are
closely related. The integrals of polar 6-complex 
functions are independent of path
in regions where the functions are regular. The fact that the exponential form
of the polar 6-complex numbers depends on the cyclic variables $\phi_1, \phi_2$
leads to the 
concept of pole and residue for integrals on closed paths. The polynomials of
polar 6-complex variables can be written as products of linear or quadratic
factors.

The operations of addition and multiplication of the planar  6-complex numbers
introduced in this 
work have a geometric interpretation based on the amplitude $\rho$,
the modulus $d$, the planar angles $\psi_1, \psi_2$ and the azimuthal angles
$\phi_1, \phi_2, \phi_3$.  
The planar 6-complex numbers can be written in exponential and
trigonometric forms with the aid of these variables.
The planar 6-complex functions defined by series of powers are analytic, and 
the partial derivatives of the components of the 6-complex functions are
closely related. The integrals of planar 
6-complex functions are independent of path
in regions where the functions are regular. The fact that the exponential form
of the 6-complex numbers depends on the cyclic variables $\phi_1, \phi_2,
\phi_3$ leads to the 
concept of pole and residue for integrals on closed paths. The polynomials of
planar 6-complex variables can always be written as products of linear factors,
although the factorization is not unique.

\newpage

FIGURE CAPTIONS\\

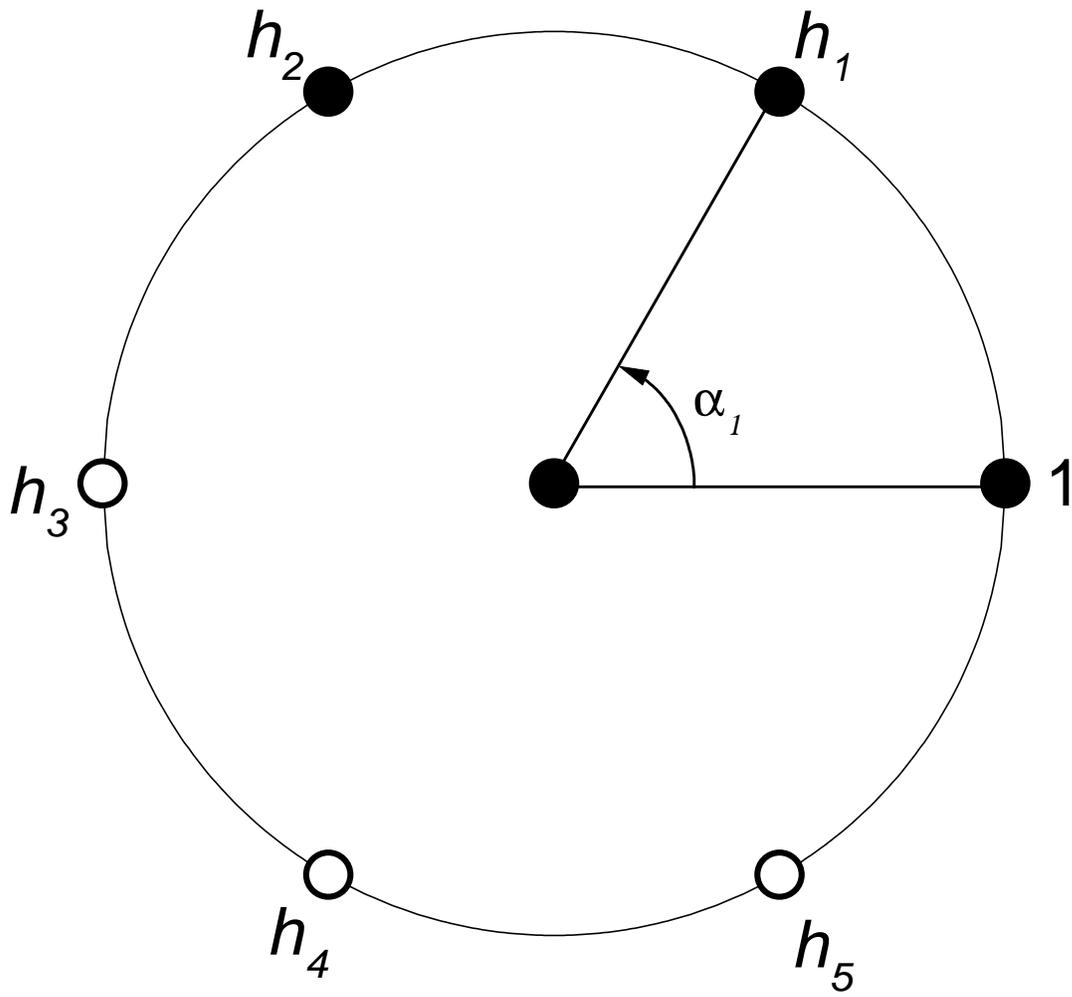
Fig. 1. Representation of the polar 
hypercomplex bases $1,h_1,h_2,h_3,h_4,h_5$
by points on a circle at the angles $\alpha_k=2\pi k/6$.
The product $h_j h_k$ will be represented by the point of the circle at the
angle $2\pi (j+k)/6$, $i,k=0,1,...,5$, where $h_0=1$. If $2\pi\leq 
2\pi (j+k)/6\leq 4\pi$, the point represents the basis
$h_l$ of angle $\alpha_l=2\pi(j+k)/6-2\pi$.\\

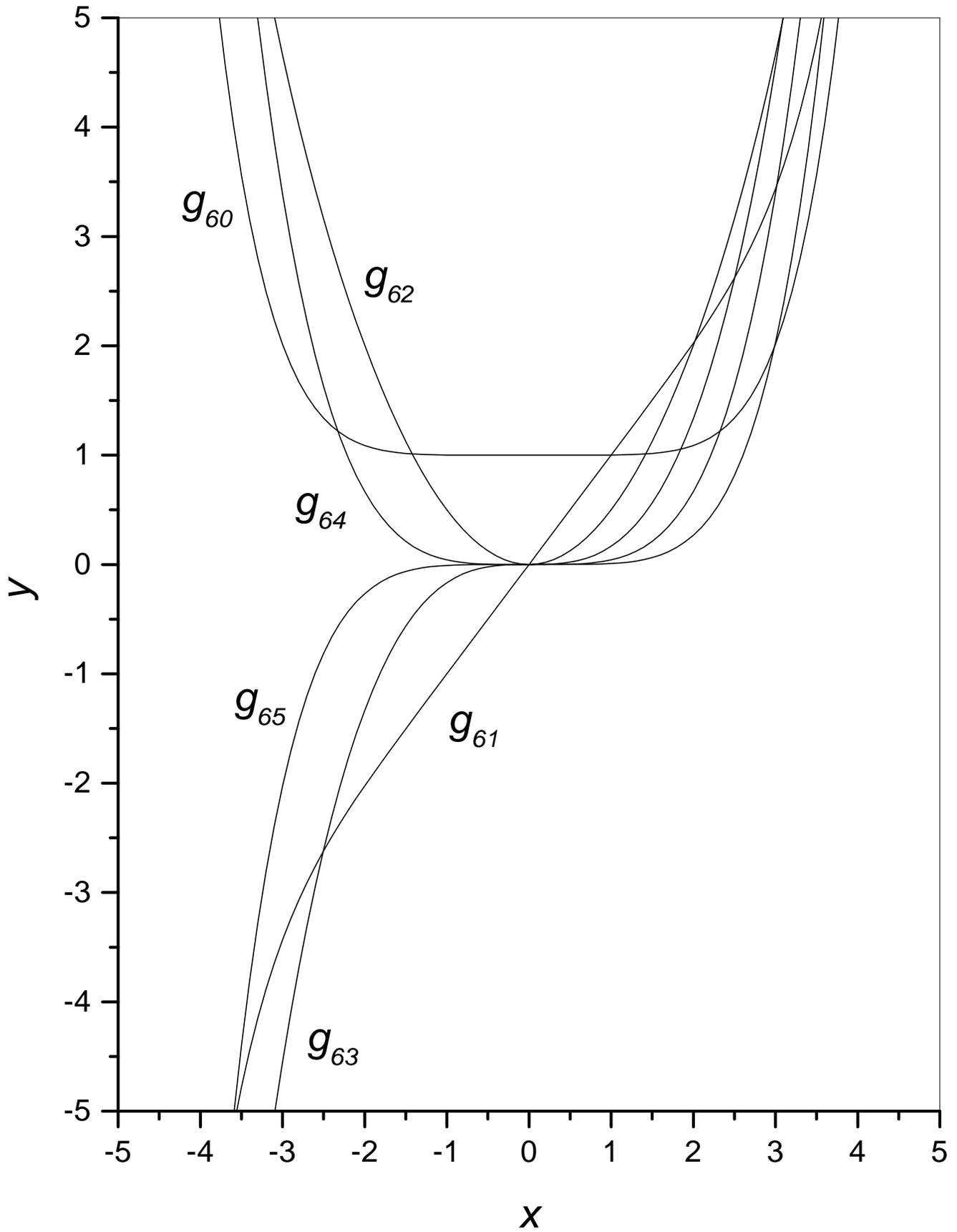
Fig. 2. Polar cosexponential functions 
$g_{60}, g_{61},g_{62}, g_{63},g_{64}, g_{65}$.\\

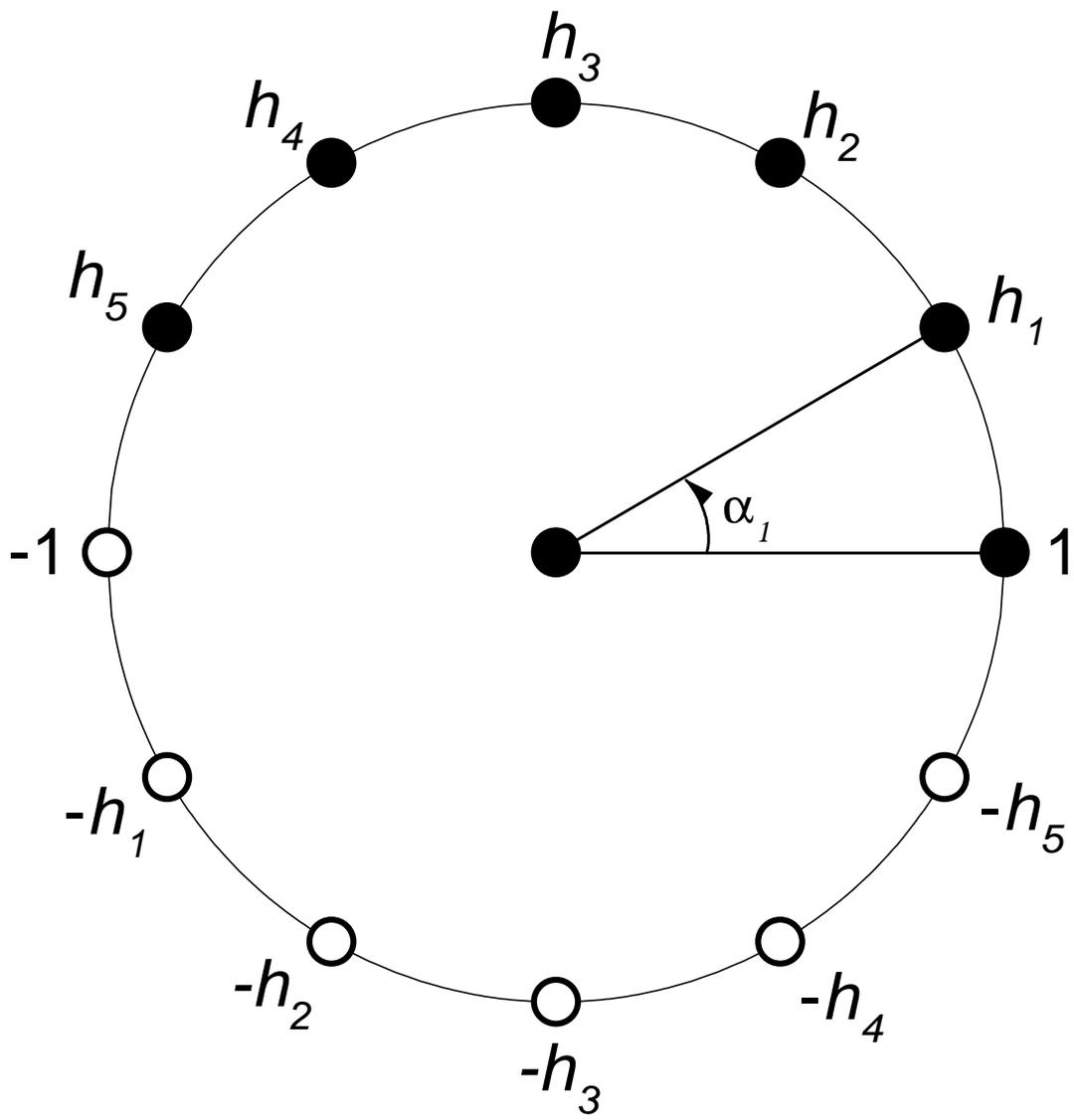
Fig. 3. Representation of the planar 
hypercomplex bases $1,h_1,h_2,h_3,h_4,h_5$
by points on a circle at the angles $\alpha_k=\pi k/6$.
The product $h_j h_k$ will be represented by the point of the circle at the
angle $\pi (j+k)/12$, $i,k=0,1,...,5$. If $\pi\leq\pi (j+k)/12\leq 2\pi$, the
point is opposite to the basis $h_l$ of angle $\alpha_l=\pi (j+k)/6-\pi$. \\

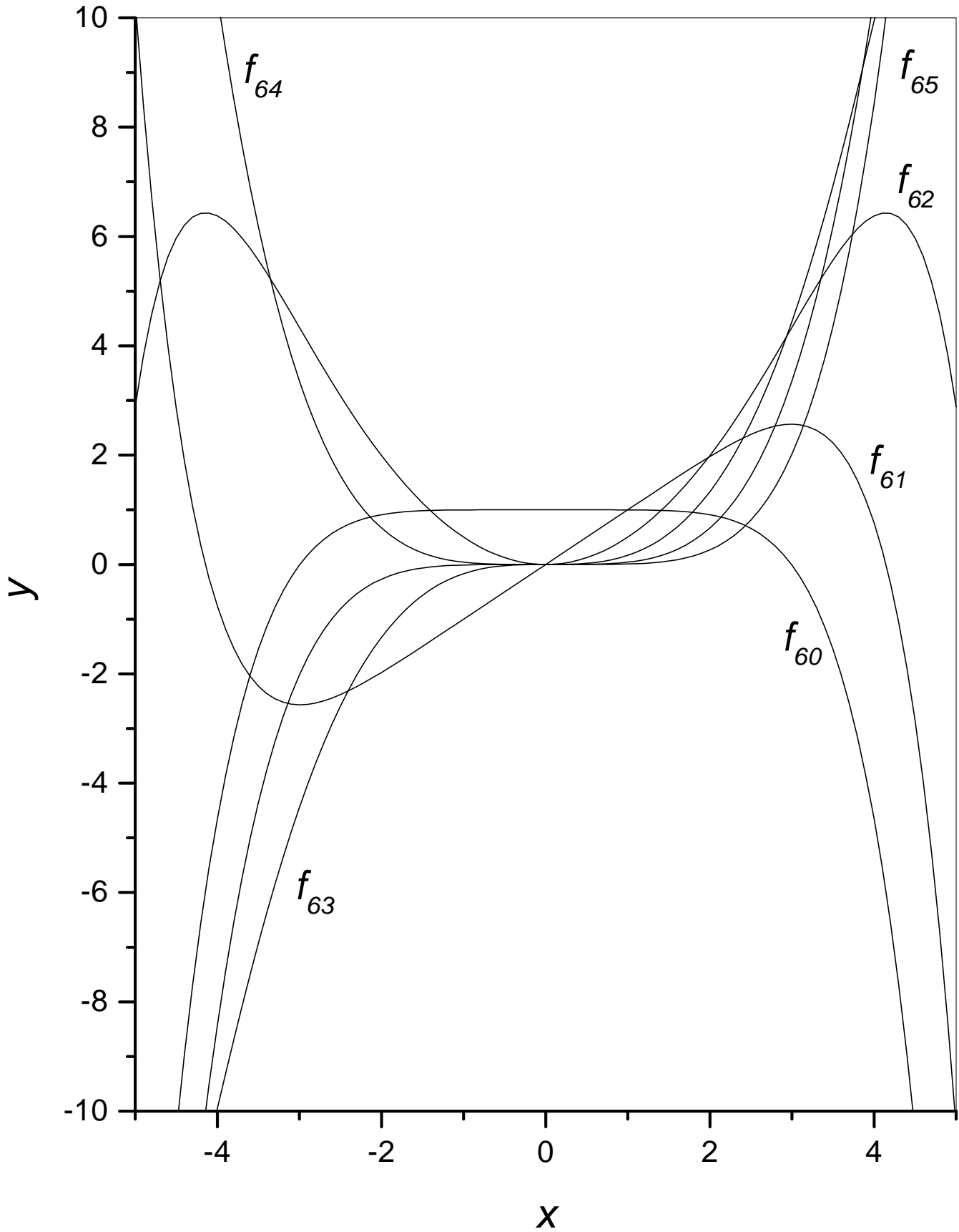
Fig. 4. Planar cosexponential functions 
$f_{60}, f_{61},f_{62}, f_{63},f_{64}, f_{65}$.\\


\begin{thebibliography}{9}
\bibitem{1} G. Birkhoff and S. MacLane, {\it Modern Algebra} (Macmillan, New
York, 
Third Edition 1965), p. 222.
\bibitem{2a} B. L. van der Waerden, {\it Modern Algebra} (F. Ungar, New York, 
Third Edition 1950),
vol. II, p. 133.
\bibitem{2} O. Taussky, Algebra, in {\it Handbook of Physics}, edited by E. U.
Condon and H. Odishaw (McGraw-Hill, New York, Second Edition 1958), p. I-22.
\bibitem{2b} 
D. Kaledin, arXiv:alg-geom/9612016;
K. Scheicher, R. F. Tichy, and K. W. Tomantschger, Anzeiger Abt.
II 134, 3 (1997); 
S. De Leo and P. Rotelli, arXiv:funct-an/9701004, 9703002;
M. Verbitsky, arXiv:alg-geom/9703016;
S. De Leo, arXiv:physics/9703033;
J. D. E. Grant and I. A. B. Strachan, arXiv:solv-int/9808019;
D. M. J. Calderbank and P. Tod, arXiv:math.DG/9911121;
L. Ornea and P. Piccinni, arXiv:math.DG/0001066.
\bibitem{2c} S. Olariu, 
{\it Hyperbolic complex numbers in two dimensions}, arXiv:math.CV/0008119;\\
{\it Complex numbers in three dimensions}, arXiv:math.CV/0008120;\\
{\it Commutative complex numbers in four dimensions}, arXiv:math.CV/0008121;\\
{\it Complex numbers in 5 dimensions}, arXiv:math.CV/0008122;\\
{\it Complex numbers in 6 dimensions}, arXiv:math.CV/0008123;\\
{\it Polar complex numbers in $n$ dimensions}, arXiv:math.CV/0008124;\\
{\it Planar complex numbers in even $n$ dimensions}, arXiv:math.CV/0008125.
\bibitem{2d} 
S. Olariu, {\it Exponential forms and path integrals for
complex numbers in $n$ dimensions}, arXiv:math.OA/0007180.
\bibitem{3} E. T. Whittaker and G. N. Watson {\it A Course of Modern
Analysis}, (Cambridge University Press, Fourth Edition 1958), p. 83.
\bibitem{4} E. Wigner, {\it Group Theory} (Academic Press, New York, 1959), p.
73.

\end{thebibliography}
\end{document}